\documentclass[10pt]{colt2018} 

\usepackage{cleveref}
\usepackage{autonum}
\DontPrintSemicolon
\usepackage{setspace,nicefrac,mathrsfs,dsfont}

\newcommand{\pushQED}[1]{%
  \toks@{\qed@elt{#1}}\@temptokena\expandafter{\QED@stack}%
  \xdef\QED@stack{\the\toks@\the\@temptokena}%
}
\newcommand{\popQED}{%
  \begingroup\let\qed@elt\popQED@elt \QED@stack\relax\relax\endgroup
}

\def\popQED@elt#1#2\relax{#1\gdef\QED@stack{#2}}
\newcommand{\qedhere}{%
  \begingroup \let\mathqed\math@qedhere
    \let\qed@elt\setQED@elt \QED@stack\relax\relax \endgroup
}

\providecommand{\proofname}{Proof}
\newtheorem{lem}{Lemma}
\newtheorem{prop}{Proposition}

\newenvironment{proofof}[1][\proofname]%
	{\par\noindent{\bfseries\upshape #1.\ }%
	}
	{\jmlrQED}

\def\hat{\widehat}
\def\tilde{\widetilde}
\newcommand{\KL}{D_{\text{\rm KL}}}

\def\RR{\mathbb R}

\def\bfI{\mathbf I}
\def\bfA{\mathbf A}
\def\bfB{\mathbf B}
\def\bfF{\mathbf F}
\def\bfW{\mathbf W}
\def\bfH{\mathbf H}

\def\bfD{\mathbf D}
\def\bfM{\mathbf M}

\def\bfY{\mathbf Y}
\def\bfU{\mathbf U}
\def\bfV{\mathbf V}
\def\bfE{\mathbf E}
\def\bfP{\mathbf P}

\def\calX{\mathcal X}

\def\calF{\mathcal F}

\def\bs{\boldsymbol}

\def\be{\bs e}
\def\bss{\bs s}
\def\bx{\bs x}
\def\by{\bs y}

\def\bu{\bs u}

\def\bY{\bs Y}

\def\bzeta{\bs\zeta}
\def\bgamma{\bs\gamma}

\def\bxi{\bs\xi}

\def\bfSigma{\bs\Sigma}

\def\tr{{\rm tr}}

\title[EWA in matrix regression and low rank]{Exponential weights in
multivariate regression and  a low-rankness favoring prior}
\usepackage{times}

 \coltauthor{\Name{Arnak S. Dalalyan} \Email{arnak.dalalyan@ensae.fr}\\
 \addr ENSAE-CREST}

\begin{document}

\maketitle

\begin{abstract}
We establish theoretical guarantees for the expected prediction error of the exponential
weighting aggregate in the case of multivariate regression that is when the label vector is
multidimensional. We consider the regression model with fixed design and bounded noise.
The first new feature uncovered by our guarantees is that it is not necessary to require
independence of the observations: a symmetry condition on the noise distribution alone
suffices to get a sharp risk bound. This result needs the regression vectors to be bounded.
A second curious finding concerns the case of unbounded regression vectors but independent
noise. It turns out that applying exponential weights to the label vectors perturbed by a
uniform noise leads to an estimator satisfying a sharp oracle inequality. The last
contribution is the instantiation of the proposed oracle inequalities to problems in which
the unknown parameter is a matrix. We propose a low-rankness favoring prior and show that
it leads to an estimator that is optimal under weak assumptions.
\end{abstract}


\section{Introduction}

The goal of this paper is to extend the scope of the applications of  the exponentially
weighted aggregate (EWA) to regression problems with multidimensional labels. Such an extension
is important since it makes it possible to cover such problems as the multitask learning, the
multiclass classification and the matrix factorization. We consider
the regression model with fixed design and additive noise. Our main contributions are
mathematical: we establish risk bounds taking the form of PAC-Bayesian type oracle
inequalities under various types of assumptions on the noise distribution.

Sharp risk bounds for the exponentially weighted aggregate in the regression with univariate labels
have been established in \citep{Leung2006,DT07,DT08,DT12b}. These bounds hold under various assumptions
on the noise distribution and cover popular examples such as Gaussian, Laplace, uniform and Rademacher
noise. One of the important specificities of the setting with multivariate labels is that noise is
multivariate as well, and one has to cope with possible correlations within its components. We provide
results that not only allow for dependence between noise components corresponding to different labels,
but also for dependence between different samples. The corresponding result, stated in \Cref{th:1},
requires, however, some symmetry of the noise distribution. To our knowledge, this is the first oracle
inequality that is sharp (\textit{i.e.}, the leading constant is equal to one) and valid under such
a general condition on the noise distribution. The remainder term in that inequality is of the order
$K/n$, where $K$ is the number of labels and $n$ is the sample size. This order of magnitude of the
remainder term is optimal, in the sense that when all the labels are equal we get the best possible
rate.

Nevertheless, one can expect that for weakly correlated labels the remainder term might be of
significantly smaller order. This is indeed the case, as shown in \Cref{th:2}, under the additional
hypothesis that the $n$ samples are independent. In the obtained sharp oracle inequality, the remainder
term is now proportional to $\|\bfSigma\|/n$, where $\|\bfSigma\|$ is the spectral norm of the noise
covariance matrix $\bfSigma\in\RR^{K\times K}$. Of course, when all the components of the noise vector
are highly correlated, the spectral norm $\|\bfSigma\|$ is proportional to $K$ and, therefore, the
conclusions of \Cref{th:1} and \Cref{th:2} are consistent.

The two aforementioned theorems are established under the condition that the aggregated matrices
belong to a set having a bounded diameter. The resulting risk bounds scale linearly in that diameter
and eventually blow up when the diameter is equal to infinity. However, it has been noticed in
that for some distributions this condition can be dropped without
deteriorating the remainder term. In  particular, this is the case of the Gaussian
\citep{Leung2006} and the uniform distributions \citep{DT08}. Furthermore, using an extended version
of Stein's lemma, \citep{DT08} show that the same holds true for any distribution having
bounded support and a density bounded away from zero. Corollary 1 in \citep{DT08} even claims that
the same type of bound holds for any symmetric distribution with bounded support. Unfortunately, the
proof of this claim is flawed since it relies on Lemma 3 (page 58) that is wrong\footnote{See \Cref{app:B} 
for more details.}. In the present work,
we have managed to repair this shortcoming and to establish sharp PAC-Bayesian risk bounds for
any symmetric distribution with bounded support. This is achieved using a key modification of the
aggregation procedure, which consists in adding a suitable defined uniform noise to data vectors
before applying the exponential weights. We call the resulting procedure noisy exponentially weighted
aggregate. Its statistical properties are presented in \Cref{th:4}.

Finally, we show an application of the obtained PAC-Bayes inequalities to the case of low-rank matrix
estimation. We exhibit a well suited prior distribution, termed spectral scaled Student prior, for
which the PAC-Bayes inequality leads to optimal remainder term. This prior is the matrix analogue of
the scaled Student prior studied in \citep{DT12b}. We also provide some hints how this estimator
can be implemented using the Langevin Monte Carlo algorithm and present some initial experimental results
on the problem of digital image denoising.

\paragraph{Notation} For every integer $k\ge 1$,
we write $\mathbf 1_k$ (resp.\ $\mathbf 0_k$) for the vector of $\RR^k$ having all
coordinates equal to one (resp.\ zero). We set $[k]=\{1,\ldots,k\}$. For every $q\in[0,\infty]$,
we denote by $\|\bu\|_q$ the usual $\ell_q$-norm of $\bu\in\RR^k$, that is
$\|\bu\|_q =(\sum_{j\in [k]}|u_j|^q)^{1/q}$ when $0<q<\infty$, $\|\bu\|_0 =
\text{Card}(\{j:u_j\not=0\})$ and $\|\bu\|_\infty =
\max_{j\in[k]} |u_j|$.

For all integers $p\ge 1$, ${\bf I}_p$ refers to the identity matrix in $\RR^{p\times p}$.
Finally the transpose and the Moore-Penrose pseudoinverse of a matrix $\bfA$ are denoted by
$\bfA^\top$ and $\bfA^\dag$, respectively. The spectral norm, the Fobenius norm and the nuclear
norm of $\bfA$ will be respectively denoted by $\|\bfA\|$, $\|\bfA\|_F$ and $\|\bfA\|_1$. For every
integer $k$, $t_k$ and $\chi^2_k$ are the Student and the chi-squared distributions with $k$
degrees of freedom.

\section{Exponential weights for multivariate regression}

In this section we describe the setting of multivariate regression and the main principles
of the aggregation by exponential weighting.

\subsection{Multivariate Regression Model}
We consider the model of multivariate regression with fixed design, in which we observe
$n$ feature-label pairs $(\bx_i,\bY_i)$, for $i\in[n]$. The labels $\bY_i\in\RR^K$
are random vectors with real entries, the features are assumed to be
deterministic elements of an arbitrary space $\calX$. Note that, unless specified otherwise,
we do not assume that the observations are independent.

We introduce the regression function $f^*:\calX\to\RR^K$ and noise vectors $\bxi_i$:
\begin{align}
\bfF_i^* = \bfE[\bY_i] = f^*(\bx_i),\qquad \bxi_i = \bY_i-\bfF_i^*,\qquad i\in[n].
\end{align}
We are interested in estimating the values of $f^*$ at the points $\bx_1,\ldots,\bx_n$ only, which
amounts to denoising the observed labels $\bY_i$. In such a setting, of course, one can
forget about the features $\bx_i$ and the function $f^*$, since the goal is merely to estimate
the $K\times n$ matrix $\bfF^*= [\bfF_1^*,\ldots,\bfF_n^*]$.
The quality of an estimator $\hat\bfF$ will be measured using the empirical loss
\begin{align}
\ell_n(\hat\bfF,\bfF^*) = \frac1{n} \|\hat\bfF-\bfF^*\|_{F}^2 = \frac1n\sum_{i\in[n]}
\|\hat\bfF_i-\bfF^*_i\|_2^2.
\end{align}
This quantity is also referred to as in-sample prediction error. The following assumption will be
repeatedly used.

\medskip

\noindent
\textbf{Assumption C$(B_\xi,L)$.}
\textit{For some positive numbers $B_\xi$ and $L$ that, unless otherwise specified, may be equal
to $+\infty$, it holds that
\begin{align}\label{cond:1}
\max_{i\in[n]}\bfP\big(\|\bxi_i\|^2_2> KB_\xi^2\big) =0,\qquad
\sup_{\bfF,\bfF'\in\calF}\max_{i\in[n]}
\|\bfF_i-\bfF'_i\|_2^2\le KL^2.
\end{align}}

\medskip
Note in \eqref{cond:1} the presence of the normalizing factor $K$ in the upper bounds
on the Euclidean norms of $K$-dimensional vectors $\bxi_i$ and $(\bfF-\bfF')_i$. This
allows us to think of the constants $B_\xi$ and $L$ as dimension independent quantities.

\begin{figure}
\hglue-15pt\includegraphics[width=1.05\textwidth,height = 0.25\textwidth]{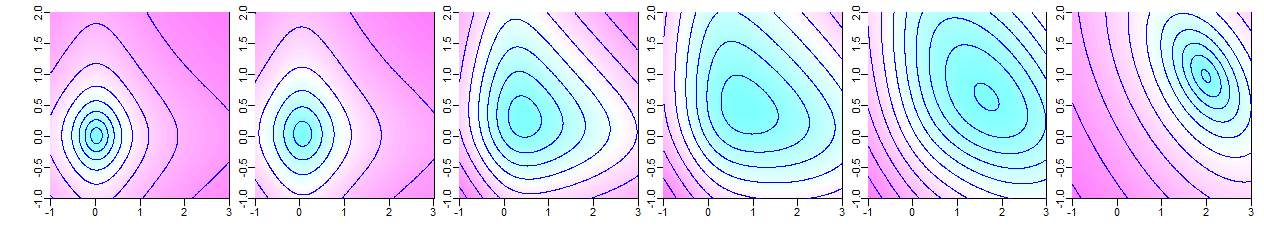}
\vglue-10pt
\caption{The contour plots of the log pseudo-posterior for different values of the
temperature parameter. The prior is a product of two Student $t(3)$ distributions.
For a very large temperature, $\tau = 20$, the first plot from the left, the posterior
is very close to the prior. On the other extreme, for $\tau = .008$, the utmost right
plot, the posterior gets close to a Dirac mass at the observed data $\bfY$ (here
$\bfY=[2,1]$).}
\vspace{-20pt}
\end{figure}

\subsection{Exponentially weighted aggregate}

The exponentially weighted aggregate (EWA) is defined as the average with respect to
a tempered posterior distribution $\pi_n$ on $\calF$, the set of all $K\times n$ matrices
with real entries. To define the tempered posterior $\pi_n$, we choose a prior distribution
$\pi_0$  on $\calF$ and a temperature parameter $\tau>0$, and set
\begin{align}
\pi_n(d\bfF) \propto \exp\Big\{-(\nicefrac1{2\tau}) \ell_n(\bfF,\bfY)\Big\}\,\pi_0(d\bfF).
\end{align}
The EWA is then
\begin{align}\label{EWA}
\hat \bfF{}^{\rm EWA} = \int_{\calF} \bfF\,\pi_n(d\bfF).
\end{align}
According to the Varadhan-Donsker variational formula, the posterior distribution
$\pi_n$ is the solution of the following optimisation problem:
\begin{align}
\pi_n\in \text{arg}\min_p \bigg\{\int_{\calF} \frac12\ell_n(\bfF,\bfY)\,p(d\bfF) + \tau \KL(p\,\|\,\pi_0)\bigg\},
\end{align}
where the inf is taken over all probability measures $p$ on $\calF$. We see that the
posterior distribution minimises a cost function which contains a term accounting for
the fidelity to the observations and a regularisation term proportional to the divergence
from the prior distribution. The larger the temperature $\tau$, the closer the posterior
$\pi_n$ is to the prior $\pi_0$.

In most situations the integral in \eqref{EWA} cannot be computed in closed form. Even its
approximate evaluation using a numerical scheme is often difficult. An appealing alternative
is then to use Monte Carlo integration. This corresponds to drawing $N$ samples $\bfF_1,
\ldots,\bfF_N$ from the posterior distribution $\pi_n$ and to define the Monte Carlo version
of the EWA by
\begin{align}
\hat\bfF^{\text{\rm MC-EWA}} = \frac1N\sum_{\ell = 1}^N \bfF_\ell.
\end{align}
Of course, the applicability of this method is limited to distributions $\pi_n$ for which
the problem of sampling can be solved at low computational cost. We will see below that
this approximation satisfies the same kind of oracle inequality as the original EWA.

\section{PAC-Bayes type risk bounds}

In this section, we state and discuss several risk bounds for the EWA and related estimators
under various conditions. We start with the case of the bounded regression vectors, \textit{i.e.},
the case where the constant $L$ in Assumption C$(B_\xi,L)$ is finite.

\subsection{Bounds without independence assumption but finite $L$}

We first state the results that hold even when the columns and rows of the noise matrix
$\bxi$ are dependent. These results, however, require the boundedness of the set of aggregated
elements $\bfF$.

\begin{theorem}\label{th:1}
Suppose that Assumption C$(B_\xi,L)$ is satisfied and the distribution
of $\bxi$ is symmetric in the sense that for any sign vector $\bss\in\{\pm 1\}^n$,
the equality in distribution
$[s_1\bxi_1,\ldots,s_n\bxi_n]\stackrel{\mathscr D}{=} \bxi$ holds. Then, for every
$\tau \ge (\nicefrac1n)(KB_\xi)(2L\vee 3B_\xi)$, we have
\begin{align}\label{GOI}
\bfE[\ell_n(\hat\bfF{}^{\rm EWA},\bfF^*)]\le
  \inf_{p}\bigg\{\int_{\calF} \ell_n(\bfF,\bfF^*)\,p(d\bfF) +
2\tau\KL(p\,||\pi_0)\bigg\},
\end{align}
where the inf is taken over all probability measures on $\calF$. Furthermore,
for larger values of the temperature, $\tau \ge (\nicefrac1n)(KB_\xi)(2L\vee 6B_\xi)$,
the following upper bound holds for $\hat\bfF = \hat\bfF{}^{\rm EWA}$
\begin{align}
\bfE[\ell_n(\hat\bfF,\bfF^*)]\le
  \inf_{p}\bigg\{\int_{\calF} \ell_n(\bfF,\bfF^*)\,p(d\bfF) +
2\tau\KL(p\,||\pi_0)\bigg\}-\frac12\int_{\calF}\bfE[\ell_n(\hat\bfF,\bfF)\pi_n(d\bfF)].
\end{align}
\end{theorem}

One can remark that the risk bound provided by \eqref{GOI} is an increasing function
of the temperature. Therefore, the best risk bound is obtained for the smallest allowed
value of temperature, that is
\begin{align}
\tau = \frac{K}{n}\,B_\xi(2L\vee 3 B_\xi).
\end{align}
Assuming $B_\xi$ and $L$ as constants, while $K = K_n$ can grow with $n$, we see that the
remainder term in \eqref{GOI} is of the order $K/n$. We will see below that using other
proof techniques, under somewhat different assumptions on the noise distribution, we can
replace $K$ by the spectral norm of the noise covariance matrix $\bfE[\bxi_i\bxi_i^\top]$.
In the ``worst case'' when all the entries of $\bxi_i$ are equal, these two bounds are
of the same order since $\|\bfE[\bxi_i\bxi_i^\top]\| = \bfE[\xi_{i1}^2]\|\mathbf 1_K\mathbf
1_K^\top\|= K\bfE[\xi_{i1}^2]$. Note, however, that the result above does not assume any
independence condition on the noise vectors $\bxi_i$.

\begin{theorem}\label{th:2}
We assume that for some $p\times p$ matrix $\bfSigma\succeq 0$, we have $\bxi =
\bfSigma^{1/2}\bar\bxi$ where $\bar\bxi$ has independent rows $\bar\bxi_{j\bullet}$ satisfying
the following boundedness and symmetry conditions:
\begin{itemize}
\item for any $(i,j)\in [n]\times[p]$, we have $\bfP(|\bar\xi_{ji}|\le \bar B_{\xi})=1$,
\item for any sign vector $\bss\in\{\pm 1\}^n$,
the equality in distribution
$[s_1\bar\bxi_{j,1},\ldots,s_n\bar\bxi_{j,n}]\stackrel{\mathscr D}{=} \bar\bxi_{j\bullet}$ holds.
\end{itemize}
In addition, the set $\calF$ is such that for some $\bar L>0$,
we have $\max_{i\in[n]}\|\bfSigma^{1/2}(\bfF_{i}-\bfF'_{i})\|_{\infty}\le \bar L$ for every
$\bfF,\bfF'\in\calF$.  Then, for every
$\tau \ge (\nicefrac1n)(\bar B_{\xi})(2\bar L\vee 3\|\bfSigma\|\bar B_{\xi})$, we have
\begin{align}\label{GOI'}
\bfE[\ell_n(\hat\bfF{}^{\rm EWA},\bfF^*)]\le
  \inf_{p}\bigg\{\int_{\calF} \ell_n(\bfF,\bfF^*)\,p(d\bfF) +
2\tau\KL(p\,||\pi_0)\bigg\},
\end{align}
where the inf is taken over all probability measures on $\calF$. Furthermore,
for larger values of the temperature, $\tau \ge (\nicefrac1n)(\bar B_{\xi})(2\bar
L\vee 6\bar B_{\xi})$, the following upper bound holds for $\hat\bfF = \hat\bfF{}^{\rm EWA}$
\begin{align}
\bfE[\ell_n(\hat\bfF,\bfF^*)]\le
  \inf_{p}\bigg\{\int_{\calF} \ell_n(\bfF,\bfF^*)\,p(d\bfF) +
2\tau\KL(p\,||\pi_0)\bigg\}-\frac12\int_{\calF}\bfE[\ell_n(\hat\bfF,\bfF)\pi_n(d\bfF)].
\end{align}
\end{theorem}

The strength of this theorem is that it does not require the independence of the observations
$\bY_i$ corresponding to different values of $i\in[n]$. Only a symmetry condition is required.
Furthermore, the resulting risk bound is valid for a temperature parameter which is of order
$O(1/n)$ and, hence, is independent of the dimension $K$ of label vectors $\bY_i$.

The proofs of \Cref{th:1} and \Cref{th:2}, postponed to \Cref{sec:proofs}, rely on the following
interesting construction related to the Skorokhod embedding. If $\gamma>0$ is a fixed number and
$\xi$ is a random variable having a symmetric distribution, then one can devise a new random variable
$\zeta$ such that $\xi+\gamma\zeta$ has the same distribution as $\xi$ and $\bfE[\zeta\,|\,\xi] =0$.
The construction of the pair $(\xi,\zeta)$ is as follows. We first draw a random variable
$R$ of the same distribution as $|\xi|$ and a Brownian motion $(B_t:t\ge 0)$ independent of $R$.
We then define the two stopping times $T$ and $T_\gamma$ by
\begin{align}
T = \inf\{t\ge 0: |B_t| = R\},\qquad T_\gamma =\inf\{t\ge 0: |B_t| = (1+\gamma)R\}.
\end{align}
One can easily check that the random variable $B_T$ has the same distribution as $\xi$ whereas
$B_{T_\gamma}$ has the same distribution as $(1+\gamma)\xi$. Furthermore, since conditionally
to $B_T=x$, the process $(B_{T+t}-x:t\ge 0)$ is a Brownian motion, we have
$\bfE[B_{T_\gamma}-B_T|B_T] = 0$. Therefore, the pair $\xi:=B_T$ and
$\zeta:=(B_{T_\gamma}-B_T)/\gamma$ satisfies the aforementioned conditions. If we set $\eta =
\zeta/\xi$, we can check that
\begin{align}
\eta =
\begin{cases}
1,&\ \text{ with probability } 1-\frac{\gamma}{1+2\gamma},\\
-1-\frac1\gamma, &\ \text{ with probability } \frac{\gamma}{1+2\gamma}.
\end{cases}
\end{align}
This is exactly the formula used in \Cref{lem:holder} below. This particular
example of the Skorokhod embedding relies heavily on the symmetry of the
distribution of $\xi$. There are other constructions that do not need this
condition. We believe that some of them can be used to further relax the
assumptions of \Cref{th:1} and \Cref{th:2}. This is, however, out of scope
of the present work.

\subsection{Bounds under independence with infinite $L$}

The previous two theorems require from the set $\calF$ of aggregated elements
to have a finite diameter $L$ (or $\bar L$) and this diameter enters (linearly)
in the risk bound through the temperature. The presence of this condition is
dictated by the techniques of the proofs; we see no reason for the established
oracle inequalities to fail in the case of infinite $L$. In the present section,
we state some results that are proved using another technique, building on the
celebrated Stein lemma, which do not need $L$ to be finite.

\begin{theorem}\label{th:3}
Assume that for some $K\times K$ positive semidefinite matrix $\bfSigma$, the
noise matrix $\bxi=\bfSigma^{1/2}\bar\bxi$ with $\bar\bxi$ satisfying the
following conditions:
\begin{itemize}\itemsep=0pt
\item[{\rm C1.}] all the random variables $\bar\xi_{j,i}$ are iid with zero mean and bounded variance,
\item[{\rm C2.}] the measure $m_{\bar\xi}(x)\,dx$, where $m_{\bar\xi}(x) = -\bfE[\bar\xi_{j,i}
\mathds 1(\bar\xi_{j,i}\le x)]$, is absolutely continuous with respect to the
distribution of $\bar\xi_{j,i}$ with a Radon-Nikodym derivative\footnote{This means
that for any bounded and measurable function $h$, we have $\int_{\RR} h(x)
m_{\bar\xi}(x)\,dx = \bfE[h(\bar\xi_{j,i})g_{\bar\xi}(\bar\xi_{j,i})]$.} $g_{\bar\xi}$,
\item[{\rm C3.}] $g_{\bar\xi}$ is bounded by some constant $G_{\bar\xi}<\infty$.
\end{itemize}
Then, for any $\tau\ge (\|\bfSigma\| G_{\bar\xi})/n$, we  have
\begin{align}\label{PAC:3}
\bfE[\ell_n(\hat\bfF{}^{\rm EWA},\bfF^*)]
		&\le \inf_p\Big\{\int_{\calF}\ell_n(\bfF,\bfF^*)\,p(d\bfF) + 2\tau \KL(p\,\|\,\pi_0)\Big\}.
\end{align}
Furthermore, if $\tau\ge 2(\|\bfSigma\| G_{\bar\xi})/n$, then
for $\hat\bfF = \hat\bfF{}^{\rm EWA}$
\begin{align}
\bfE[\ell_n(\hat\bfF,\bfF^*)]
		&=\inf_p\Big\{\int_{\calF}\ell_n(\bfF,\bfF^*)\,p(d\bfF) + 2\tau \KL(p\,\|\,\pi_0)\Big\}-
			\frac12\int_{\calF}
			\bfE\big[\ell_n(\bfF,\hat\bfF)\,\pi_n(d\bfF)\big].
\end{align}
\end{theorem}

As mentioned in \citep[pp 43-44]{DT08}, many distributions satisfy assumptions
C2 and C3. For instance, for the Gaussian distribution $\mathcal N(0,\sigma^2)$
and for the uniform in $[-b,b]$ distributions these assumptions are fulfilled
with $G_{\bar\xi} = \sigma^2$ and $G_{\bar\xi} = b^2/2$, respectively. More generally,
if $\bar\xi_{j,i}$ has a density $p_{\bar\xi}$ with bounded support $[-b,b]$, then
the assumptions are satisfied with $G_{\bar\xi} = \bfE[|\bar\xi_{j,i}|]/
\min_{|x|\le b} p_{\bar\xi}(x)$. Here, we add another class to
the family of distributions satisfying C2 and C3: unimodal distributions with
compact support.

\begin{prop}\label{prop:1}
Assume that $\bar\xi_{j,i}$ has a density $p_{\bar\xi}$ with respect to the Lebesgue measure
such that $p_{\bar\xi}(x) = 0$ for every $x\not\in[-b,b]$ and, for some $a\in [-b,b]$,
$p_{\bar\xi}$ is increasing on $[-b,a]$ and decreasing on $[a,b]$. Then, $\bar\xi_{j,i}$
satisfies C2 and C3 with $G_{\bar\xi} = (\nicefrac{b^2}{2})$.
\end{prop}

Perhaps the most important shortcoming of the last theorem is that it cannot be applied to
the discrete distributions of noise. In fact, if the distribution of $\bar{\xi}_{j,i}$
is discrete, then there is no chance condition C2 to be satisfied. This is due to the fact that
the measure $m_{\bar{\xi}}\,dx$, being absolutely continuous with respect to the Lebesgue
measure, cannot be absolutely continuous with respect to a counting measure. On the other hand,
\Cref{th:1} and \Cref{th:2} can be applied to discrete noise distributions, but they require
boundedness of the family $\calF$. At this stage, we do not know whether it is possible to
extend PAC-Bayesian type risk bound \eqref{PAC:3} to discrete distributions
and unbounded sets $\calF$. However, in the case of a bounded discrete noise, we propose
a simple modification of the EWA for which \eqref{PAC:3} is valid.

The modification mentioned in the previous paragraph consists in adding a uniform noise
to the entries of the observed labels $\bfY_i$. Thus, we define the noisy exponential weighting
aggregate, nEWA, by
\begin{align}\label{nEWA}
\hat\bfF{}^{\rm nEWA} = \int_{\calF} \bfF\,\bar\pi_n(d\bfF),\qquad
\bar{\pi}_n(d\bfF)\propto \exp\big\{-(\nicefrac1{2\tau})\ell_n(\bfF,\bar{\bfY})\big\}\,\pi_0(d\bfF),
\end{align}
where $\bar{\pi}_n$ is defined in the same way as $\pi_n$ but for the perturbed matrix
$\bar{\bfY} = \bfY + \bzeta$, with $\bzeta$ a $K\times n$ random perturbation matrix.

\begin{theorem}
	\label{th:4}
Let $\hat\bfF{}^{\rm nEWA}$ be the noisy EWA defined by \eqref{nEWA}.
Assume that
\begin{itemize}\itemsep=0pt
	\item[{\rm C4.}] entries $\xi_{j,i}$ of the noise matrix $\bxi$ are iid with zero mean and bounded by
	some constant $B_{\xi}>0$,
	\item[{\rm C5.}] entries $\zeta_{j,i}$ of the perturbation matrix are iid uniformly distributed in
	$[-B_{\xi},B_{\xi}]$.
\end{itemize}
Then, for any $\tau\ge 2B_\xi^2/n$, we  have
\begin{align}\label{PAC:4}
\bfE[\ell_n(\hat\bfF{}^{\rm nEWA},\bfF^*)]
&\le \inf_p\Big\{\int_{\calF}\ell_n(\bfF,\bfF^*)\,p(d\bfF) + 2\tau \KL(p\,\|\,\pi_0)\Big\}.
\end{align}
Furthermore, if $\tau\ge 4B_\xi^2/n$, then
for $\hat\bfF = \hat\bfF{}^{\rm nEWA}$
\begin{align}
\bfE[\ell_n(\hat\bfF,\bfF^*)]
&=\inf_p\Big\{\int_{\calF}\ell_n(\bfF,\bfF^*)\,p(d\bfF) + 2\tau \KL(p\,\|\,\pi_0)\Big\} -
\frac12\int_{\calF}
\bfE\big[\ell_n(\bfF,\hat\bfF{}^{\rm EWA})\,\bar\pi_n(d\bfF)\big].
\end{align}
\end{theorem}

\begin{proofof}[Proof of \Cref{th:4}]
Let us check that the matrix of
perturbed labels $\bar\bfY$ satisfies the conditions of \Cref{th:3} with $\bfSigma=\bfI_K$.
To this end, we set $\tilde\xi_{j,i} = \xi_{j,i}+\zeta_{j,i}$. We will check that the distribution
of $\tilde{\xi}_{j,i}$ satisfies conditions C2 and C3 (condition C1 is straightforward).
Since the distribution of $\tilde\xi_{j,i}$ is the convolution of that of $\xi_{j,i}$
and a uniform distribution, it admits a density with respect to the Lebesgue measure which is
given by
\begin{align}
	\tilde p(x) = \frac1{2B_{\xi}}\,\bfP(|\xi_{j,i}-x|\le B_\xi)=
	\frac1{2B_{\xi}}\,\bfP(\xi_{j,i}\in[x-B_\xi, x+B_\xi]\cap [-B_\xi,B_\xi]).
\end{align}
The set  $A_x:=[x-B_\xi, x+B_\xi]\cap [-B_\xi,B_\xi]$ is empty if $|x|>2B_\xi$,
increasing on the interval
$x\in[-2B_\xi,0]$ and decreasing on the interval $x\in[0,2B_\xi]$. This implies that
the density $\tilde p$ is zero outside the interval $[-2B_\xi,2B_\xi]$ and unimodal
in this interval. Therefore, it satisfies \Cref{prop:1} with $b=2B_\xi$ and $a= 0$.
This implies that conditions C2 and C3 are fulfilled with $G_{\bar{\xi}} = 2B_\xi^2$
and $\|\bfSigma\|=1$. Thus, the conclusion of \Cref{th:3} applies and yields the claims
of \Cref{th:4}.
\end{proofof}

We can replace in \Cref{th:3} the condition C4 by $\bxi= \bfSigma^{1/2}\bar\bxi$, where
$\bar\xi_{j,i}$ are iid and bounded. In this case, the contamination added to the labels
is of the form $\bfSigma^{1/2}\bar\bzeta$, where $\bar\zeta_{j,i}$'s are iid uniform. The
claims of \Cref{th:4} remain valid, but they are of limited interest, since it is
not likely to find a situation in which the matrix $\bfSigma$ is known.

\subsection{Risk bounds for the Monte Carlo EWA}

The four theorems of the previous sections contain all two risk bounds. The first bound is, in
each case, more elegant and valid for a smaller value of the temperature than the second bound.
However, the latter appears to be more useful for getting guarantees for the Monte Carlo version
of the EWA. This is due to the fact that the additional term in the second risk bounds is
proportional to the difference of the risks between the MC-EWA and the EWA.

\begin{proposition}
If $\hat\bfF^{\text{\rm MC-EWA}}$ is the MC-EWA with $N$ Monte Carlo samples, then
\begin{align}
\bfE[\ell_n(\hat\bfF^{\text{\rm MC-EWA}},\bfF^*)]  = \bfE[\ell_n(\hat\bfF^{\rm EWA},\bfF^*)]
+\frac1N\int_{\calF}\bfE[\ell_n(\bfF,\hat\bfF^{\rm EWA})\,\pi_n(d\bfF)].
\end{align}
Therefore, if the conditions of one of the four foregoing theorems are satisfied and $\tau$
is chosen accordingly then, as soon as $N\ge 2$,
\begin{align}
\bfE[\ell_n(\hat\bfF^{\text{\rm MC-EWA}},\bfF^*)]  \le \inf_p\Big\{\int_{\calF}\ell_n(\bfF,\bfF^*)
\,p(d\bfF) + 2\tau \KL(p\,\|\,\pi_0)\Big\}.
\end{align}
\end{proposition}

The proof of this result is straightforward and, therefore, is omitted. Note that
this result bounds only the expected error, where the expectation is taken with respect
to both the noise matrix $\bxi$ and the Monte Carlo sample. Using standard concentration
inequalities, this bound can be complemented by an evaluation of the deviation between
the Monte Carlo average $\hat\bfF^{\text{\rm MC-EWA}}$ and its ``expectation''
$\hat\bfF^{\text{\rm EWA}}$. 

\subsection{Relation to previous work}

To the best of our knowledge, the first result in the spirit of the oracle inequalities
presented in foregoing sections has been established by \cite{Leung2006}, using a
technique heavily based on Stein's unbiased risk estimate for regression with Gaussian
noise developed in \citep{George86a,George86b}. The first extensions to more general noise
distributions were presented in \citep{DT07,DT08} and later on refined in \citep{DT12b}.
In all these papers, only the problem of aggregating ``frozen'' (that is independent of
the data used for the aggregation) estimators. In his PhD thesis, \cite{Leung}
proved that analogous oracle bounds hold for the problem of aggregation of shrinkage
estimators. The case of linear estimators has been explored by \cite{DS12,DRXZ,bellec2018}.
In the context of sparsity, statistical properties of exponential weights were
studied in \cite{Alquier3,RT11,RT12}.

There is also extensive literature on the exponential weights for problems with iid
obsrvations, such as the density model, the regression with random design, etc. We refer
the interested reader to \citep{Yang00, Yang00a, Catoni1, JRT,Audibert,DT12a} and the
references therein. It is useful to note here that the proof techniques used in the iid
setting and in the setting with deterministic design considered in the present work are
very different. Furthermore, the version exponential of the exponential weights used in
the iid setting involves an additional averaging step and is therefore referred to as
progressive mixture or mirror averaging.

\section{EWA with low-rank favoring priors}

To give a concrete example of application of the results established in previous section,
let us consider the so called reduced rank regression model. An asymptotic analysis of this
model goes back to \citep{IZENMAN1975248}, whereas more recent results can be found in
\citep{bunea2011, bunea2012} and the references therein.
It corresponds to assuming that
the matrix $\bfF^* = \bfE[\bfY]$ has a small rank, as compared to its maximal possible value
$K\wedge n$. Equivalently, this means that the observed $K$ dimensional vectors $\bY_1$,
\ldots,$\bY_n$ belong, up to a noise term, to a low dimensional subspace. Such problems arise,
for instance, in subspace clustering or in multi-index problems. Of course, one can estimate
the matrix $\bfF^*$ by the PCA, but it requires rather precise knowledge of the rank.

In order to get an estimator that takes advantage of the (nearly) low-rank property
of the matrix $\bfF^*$, we suggest to use the following prior
\begin{align}\label{SSSp}
\pi_0(d\bfF) \propto \text{det}(\lambda^2\bfI_K + \bfF\bfF^\top)^{-(n+K+2)/2}\,d\bfF,
\end{align}
where $\lambda>0$ is a tuning parameter. From now on, with a slight abuse of notation, we
will denote by $\pi_0(\bfF)$ the probability density function of the measure $\pi_0(d\bfF)$.
The same will be done for $p(d\bfF)$ and $\pi_n(d\bfF)$. We will refer to $\pi_0$  as the
spectral scaled Student prior, since one easily checks that
\begin{align}
\pi_0(\bfF) \propto \prod_{j=1}^K(\lambda^2 + s_j(\bfF)^2)^{-(n+K+2)/2},
\end{align}
where $s_j(\bfF)$ denotes the $j$th largest singular value of $\bfF$. We can recognize in
the last display the density function of the scaled Student $t$ evaluated at $s_j(\bfF)$.
Thus, the scaled spectral Student prior operates on the spectrum of $\bfF$ as the
sparsity favoring prior introduced in \citep{DT12b} on the vectors. Another interesting
property of this prior, is that if $\bfF\sim \pi_0$, then the marginal distributions of
the columns of $\bfF$ are scaled multivariate Studtent $t_3$.

\begin{lem}\label{lem:t3}
If\/ $\bfF$ is a random matrix having as density the function $\pi_0$, then the random vectors
$\bfF_i$ are all drawn from the $K$-variate scaled Student distribution $(\lambda/\sqrt{3})t_{3,K}$.
As a consequence, we have $\int_{\calF} \|\bfF_i\|_2^2\,\pi_0(\bfF)\,d\bfF = 	\lambda^2 K$.
\end{lem}

From a mathematical point of view, the nice feature of the aforementioned prior is that
the Kullback-Leibler divergence between $\pi_0$ and its shifted by a matrix $\bar\bfF$
version grows proportionally to the rank of $\bar\bfF$, when all the other parameters
remain fixed. This is formalized in the next result.

\begin{lem}\label{lem:KL}
Let $\bar p$ be the probability density function obtained from the prior $\pi_0$ by a
translation, $\bar p(\bfF) = \pi_0(\bfF-\bar\bfF)$. Then, for any matrix $\bar\bfF$ of
at most rank $r$, we have
\begin{align}
\KL(\bar p\,\|\,\pi_0)
		& \le 2r(n+K+2) \log\bigg(1 + \frac{\|\bar\bfF\|_F}{\sqrt{2r}\lambda}\bigg)
		\le 2r(n+K+2) \log\big(1 + {\|\bar\bfF\|}/{\lambda}\big).
\end{align}
\end{lem}

The proof of this result is deferred to the appendix. Applying this lemma, in conjunction with
\Cref{th:3}, we get a risk bound in the reduced rank regression problem which illustrates
the power of the exponential weights.

\begin{theorem}\label{th:5}
Let the noise matrix $\bxi$ and the artificial perturbation matrix $\bzeta$ satisfy the
assumptions of \Cref{th:4}. Let $\pi_0$ be the scaled spectral Student prior \eqref{SSSp}
with some tuning parameter $\lambda>0$. Then, for every $\tau\ge 2B_{\xi}^2/n$, we have
\begin{align}\label{PAC:5}
\bfE[\ell_n(\hat\bfF{}^{\rm nEWA},\bfF^*)]
&\le\inf_{\bar\bfF}\bigg\{\ell_n(\bar\bfF,\bfF^*)
+ 4r(\bar\bfF)(n+K+2)\tau \log\bigg(1 + \frac{\|\bar\bfF\|_F}{\sqrt{2r}\lambda}\bigg)\bigg\}
+K\lambda^2,
\end{align}
where $r(\bfF) =  \text{rank}(\bfF)$ and the inf is taken over all $K\times n$ matrices $\bar\bfF$.
\end{theorem}

\begin{proofof}[Proof of \Cref{th:5}]
Let us fix an arbitrary matrix $\bar\bfF$ and denote its rank by $r$. We apply \Cref{th:4}
and upper bound the inf with respect to all probability distributions $p$ by the right hand
side of \eqref{PAC:4} evaluated at the distribution $\bar p$ defined in \Cref{lem:KL}. This
yields
\begin{align}\label{PAC:5a}
\bfE[\ell_n(\hat\bfF{}^{\rm nEWA},\bfF^*)]
&\le \int_{\calF}\ell_n(\bfF,\bfF^*)\pi_0(\bfF-\bar\bfF)\,d\bfF
+ 4 r(\bar\bfF)(n+K+2)\tau \log\bigg(1 + \frac{\|\bar\bfF\|_F}{\sqrt{2r}\lambda}\bigg).
\end{align}
Using the translation invariance of the Lebesgue measure and the fact that
$\int \bfF\,\pi_0(\bfF)\,d\bfF = \mathbf 0$, we get
\begin{align}
\int_{\calF}\ell_n(\bfF,\bfF^*)\pi_0(\bfF-\bar\bfF)\,d\bfF
		= \ell_n(\bar\bfF,\bfF^*) + \frac1n\int_{\calF}\|\bfF\|_F^2\,\pi_0(\bfF)\,d\bfF.
\end{align}
Let us focus on the evaluation of the last integral. The claimed inequality follows from
the last display by applying \Cref{lem:t3} and the fact that $\|\bfF\|_F^2 = \sum_{i\in[n]}
\|\bfF_i\|_2^2$.
\end{proofof}

There are many papers using Bayesian approach to the problem of prediction with low rank
matrices, see \citep{AlquierALT,bouwmans:hal-01373013} and the references therein. All the
methods we are aware of define a prior on $\bfF$ using the following scheme: first choose a prior
on the space of triplets $(\bfU,\bfV,\bgamma)$, where $\bfU$ and $\bfV$ are orthogonal matrices and 
$\bgamma$ is a vector with nonnegative entries, and then define $\pi_0$ as the distribution
of $\bfU\bfD_{\bgamma}^2\bfV^\top$ (see, for instance, \citep{mai2015,Alquier2017}). Similar type of 
priors have been also used in the problem of tensor decomposition and completion \cite{RaiICML} but, 
to date, their statistical accuracy has not been studied.

To our knowledge, \citep{Yang2017} is the only work dealing with prior
\eqref{SSSp} in a context related to low rank matrix estimation and completion. It
proposes variational approximations to the Bayes estimator and demonstrates their
good performance on various data sets. In a sense, \Cref{th:5} provides theoretical
justification for the empirically observed good statistical properties of the prior
defined in \eqref{SSSp}.

Let us briefly discuss the inequality of \Cref{th:5}. Assume that we choose
$\tau = 2B_\xi^2/n$ and $\lambda^2 = B_\xi^2(n+K)/K$. Then, we see that \eqref{PAC:5}
handles optimally mis-specification, since it is an oracle inequality with a
leading constant 1, and the remainder term is of optimal order $r(n+K)/n$, up to a
logarithmic factor. Other oracle inequalities with nearly optimal remainder terms 
in the context of low-rank matrix estimation and completion are exposed in
\citep{Mai15,AlquierALT,Alquier2017}. However, those results are not sharp oracle 
inequalities since the factor in front of the leading term in the upper bound is
larger than 1. 

Using the properties of the scaled Student prior exposed in \Cref{lem:t3} and
\Cref{lem:KL}, one can establish oracle inequalities in other statistical problems
in which the unknown parameter is a matrix, such as matrix completion, trace regression
or multiclass classification, see  \citep{Srebro2005, Rohde11, KLT11,
CT10, CP11, BYW,GL11,NW11,NW12,Klopp14, DGP16}. This is left to future work.

\section{Implementation and a few numerical experiments}

In this section, we report the results of some proof of concept numerical
experiments. We propose to compute an approximation of the EWA with the
scaled multivariate Student prior by a suitable version of the Langevin Monte
Carlo algorithm.  To describe the letter, let us first remark that
\begin{align}\label{logpi}
\log \pi_n(\bfF) = -\frac1{2\tau} \ell_n(\bfF,\bfY) - \frac{(n+K+2)}{2}
\log \text{det}(\lambda^2\bfI_K + \bfF\bfF^\top).
\end{align}
From this relation, we can infer that
\begin{align}
-\nabla \log \pi_n(\bfF) = \frac1{n\tau} (\bfF-\bfY) + (n+K+2)
(\lambda^2\bfI_K + \bfF\bfF^\top)^{-1}\bfF.
\end{align}
The (constant-step) Langevin MC is defined by choosing an initial matrix
$\bfF_0$ and then by using the recursion
\begin{align}
\bfF_{k+1} = \bfF_{k} + h\nabla \log \pi_n(\bfF) +\sqrt{2h}\,\bfW_k,\qquad
k=0,1,\ldots,
\end{align}
where $h>0$ is the step-size and $\bfW_0,\bfW_1,\ldots$ are independent Gaussian
random matrices with iid standard Gaussian entries. For (strongly) log-concave
densities $\pi$, nonasymptotic guarantees for the LMC have been recently
established in \cite{Dal14,Durmus16}, but they do not carry over the present
case since  the right hand side of \eqref{logpi} is not concave. Our numerical
experiments show that (despite the absence of theoretical guarantees) the LMC
converges and leads to relevant results.

Note that a direct application of the Langevin MC algorithm involves
a $K\times K$ matrix inversion at each iteration. This might be costly and
can slow down significantly the algorithm. We suggest to replace this matrix
inversion by a few steps of gradient descent for a suitably chosen optimization
problem. Indeed, one can easily check that the matrix
$\bfM = (\lambda^2\bfI_K + \bfF\bfF^\top)^{-1}\bfF$ is nothing else but the
solution to the convex optimization problem
\begin{align}
\min \big\{\|\bfI_n-\bfF^\top \bfM\|_F^2 + 2\lambda^2\|\bfM\|_F^2\big\}.
\end{align}
We use ten gradient descent steps for approximating the solution of this optimization
problem. This does not require neither matrix inversion nor svd or other costly
operation. Theoretical assessment of the Langevin MC with inaccurate gradient
evaluations can be found in \cite{DalKar17}.

We applied this algorithm to the problem of image denoising. We chose an RGB
image of resolution $120\times 160$ and applied to it an additive Gaussian white
noise of standard deviation $\sigma\in\{10,30,50\}$. In order to make use of the
denoising algorithm based on the aforementioned Langevin MC, we transformed the
noisy image into a matrix of size $192\times 300$. Each row of this transformed matrix
corresponds to a patch  of size $10\times 10\times 3$ of the noisy image. The patches
are chosen to be non-overlapping in order to get a reasonable dimensionality. We expect
the result to be better for overlapping patches, but the computational cost will also
be high. The parameters were chosen as follows:
\begin{align}
\tau = 2\sigma^2/n;\qquad \lambda = 10\sigma\sqrt{(n+K)/K};
\quad h = 10;\qquad k_{\max} = 4000.
\end{align}
Note that the values of $\tau$ and $\lambda$ are suggested by our theoretical results,
while the step-size $h$ and the number of iterations of the LMC, $k_{\max}$, were chosen
experimentally. The LMC after $k_{\max}$ iterations provides one sample that is approximately drawn from the pseudo-posterior $\pi_n$. We did $N = 400$ repetitions and
averaged the obtained results for approximating the posterior mean.

\begin{figure}
{\hspace{48pt}True image \hspace{85pt} Noisy image \hspace{74pt} Restaured image}\\
\centerline{\includegraphics[width = 0.9\textwidth]{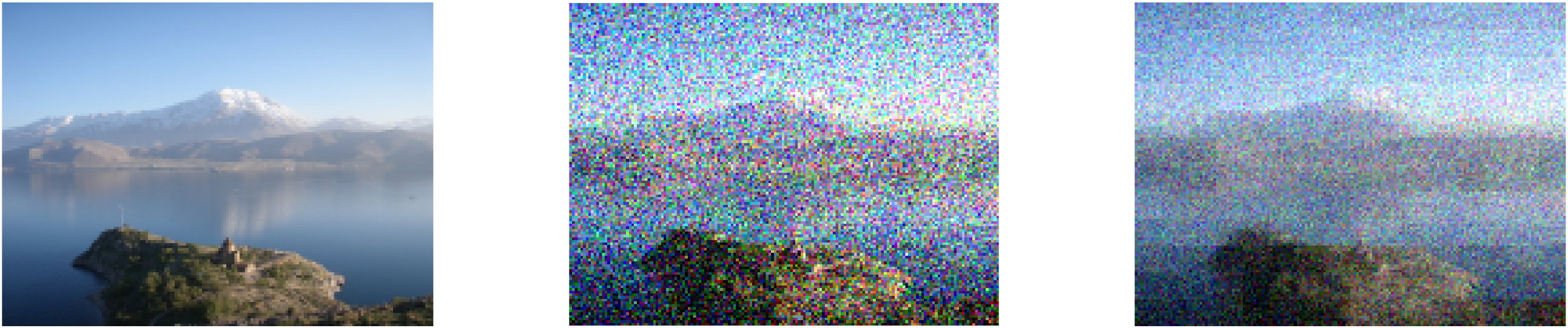}}\\
{\hphantom{True image} \hspace{114pt} $\sigma = 50$, PSNR = 14.1 \hspace{60pt} PSNR = 20.6}\\[5pt]
\centerline{\includegraphics[width = 0.9\textwidth]{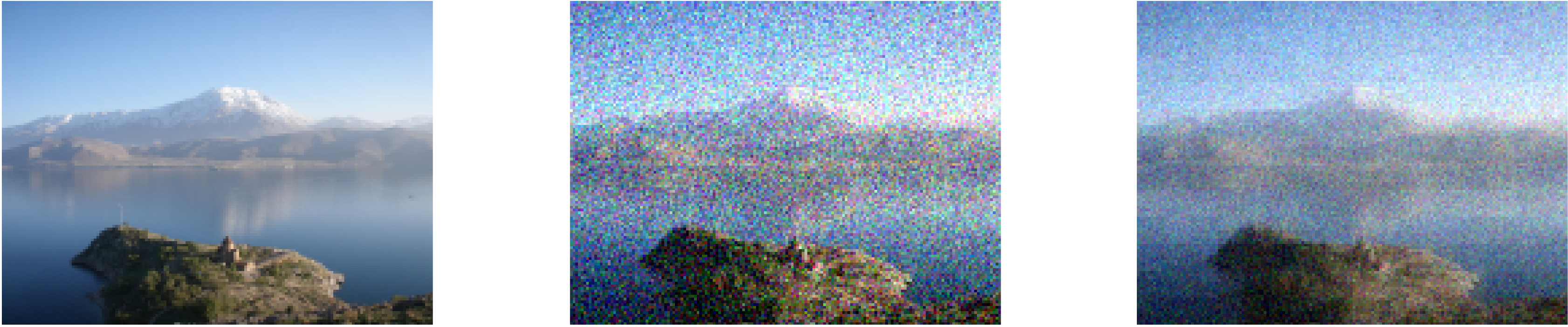}}\\
{\hphantom{True image} \hspace{114pt} $\sigma = 30$, PSNR = 18.6 \hspace{60pt} PSNR = 24.3}\\[5pt]
\centerline{\includegraphics[width = 0.9\textwidth]{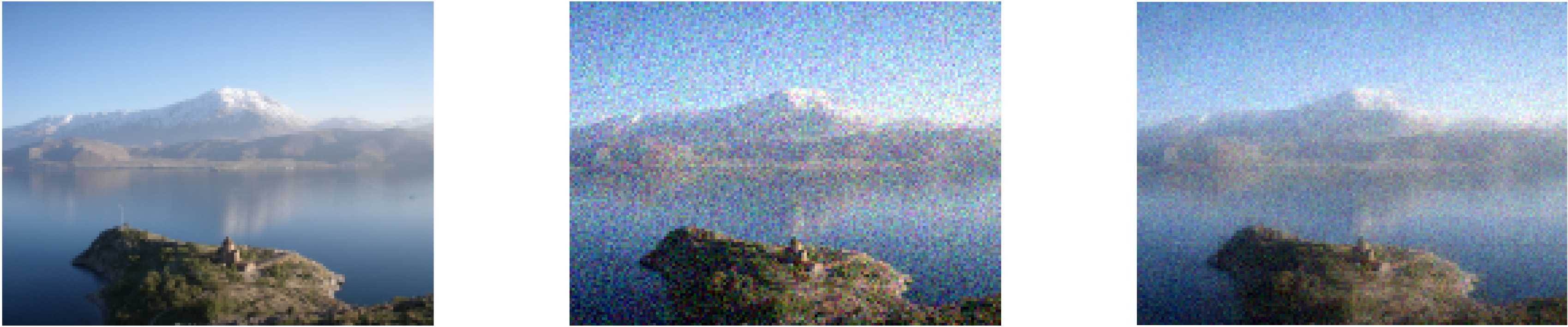}}
{\hphantom{True image} \hspace{114pt} $\sigma = 20$, PSNR = 22.1 \hspace{60pt} PSNR = 27.1}\\[5pt]
\centerline{\includegraphics[width = 0.9\textwidth]{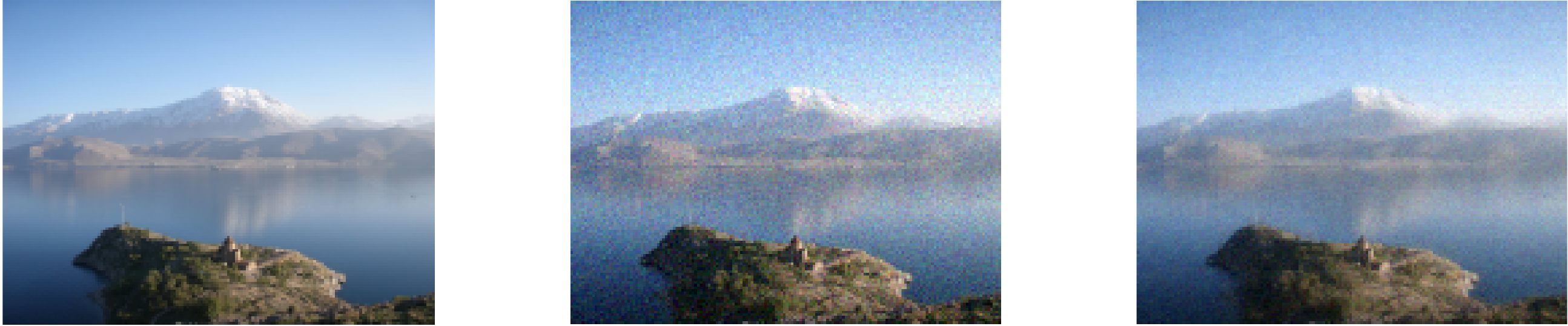}}
{\hphantom{True image} \hspace{114pt} $\sigma = 10$, PSNR = 28.1 \hspace{60pt} PSNR = 32.2}\\[5pt]
\centerline{\includegraphics[width = 0.9\textwidth]{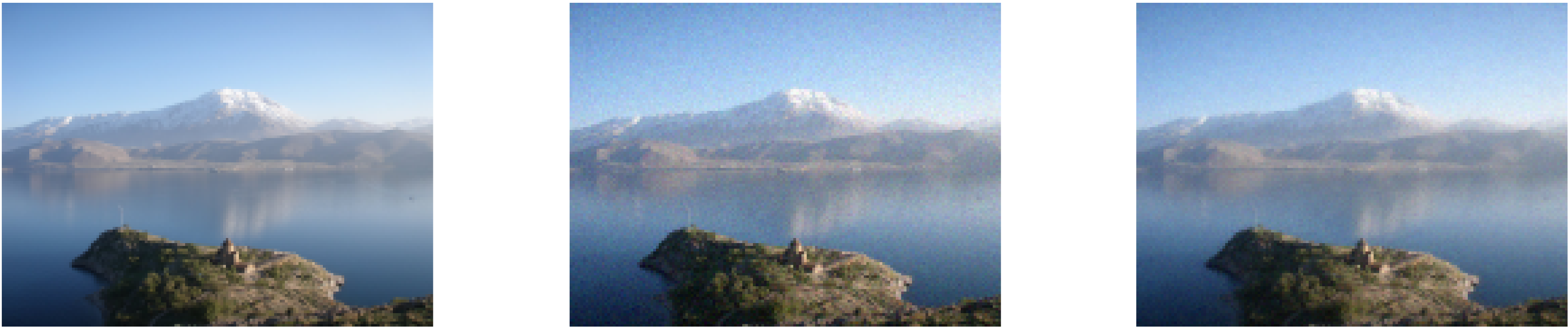}}
{\hphantom{True image} \hspace{114pt} $\sigma = 5$, PSNR = 34.1 \hspace{60pt} PSNR = 36.7}\\[5pt]
\caption{The result of the experiment on image densoising. Left: the original $120\times 160\times 3$ image.
Middle: the noisy image for different values of $\sigma$. Right: the denoised image.}
\label{fig2}
\end{figure}

%
%
%
%
%

\section{Conclusion}

We have studied the expected in-sample prediction error of the Exponentially
Weighted Algorithm (EWA) in the context of multivariate regression with possible
dependent noise. We have shown that under boundedness assumptions on the noise
and the aggregated elements, the EWA satisfies a PAC-Bayes type sharp oracle
inequality, provided that the temperature parameter is sufficiently large. The
remainder term in these oracle inequalities is of arguably optimal order of
magnitude and is consistent with the corresponding results obtained in the
model of univariate regression. An interesting observation is that if we apply
the EWA to the data matrix artificially contaminated by a uniform noise, the
resulting procedure satisfies a sharp oracle inequality under a much weaker
assumption on the noise distribution. In particular, this allows to cover any
distribution with bounded support. We have also included the results of a small
numerical experiment on image denoising, that shows the applicability of the
EWA.

\section{Proofs of the main results}\label{sec:proofs}

The proofs of all the main theorems stated in the previous sections are gathered
in this section. The proofs of some technical lemmas are deferred to \Cref{sec:lemmas}.

\begin{proofof}[Proof of \Cref{th:1}]
We wish to upper bound  $\ell_n(\hat\bfF^{\rm EWA},\bfF^*)$.
Let $\bzeta$ be a random matrix such that
$\bfE[\bzeta|\bfY] = \mathbf 0$ and define
\begin{align}
\ell_n(\bfF,\bfF^*,\bzeta)
		= \ell_n(\bfF,\bfF^*) + \frac{2}{n}\langle\bzeta,\bfF-\bfF^*\rangle.
\end{align}
In what follows, we use the short notation $\hat\bfF$ instead of $\hat\bfF^{\rm EWA}$.
We have, for every $\alpha>0$,
\begin{align}
\ell_n(\hat\bfF,\bfF^*,\bzeta)
		& = \frac1{\alpha}\log \exp\Big\{\alpha\,\ell_n(\hat\bfF,\bfF^*,\bzeta)\Big\}\\
		& = \frac1{\alpha}\underbrace{\log \int_{\calF}e^{\alpha\big(
			\ell_n(\hat\bfF,\bfF^*,\bzeta)-\ell_n(\bfF,\bfF^*,\bzeta)\big)}\,\pi_n(d\bfF)}_{:=S_1(\alpha)}
			-\frac1{\alpha}\underbrace{\log \int_{\calF} e^{-\alpha\,\ell_n(\bfF,\bfF^*,\bzeta)}\,
			\pi_n(d\bfF)
			}_{:=S(\alpha)}.
\end{align}
The next two lemmas provide suitable upper bounds on the magnitude of the terms
$S(\alpha)$ and $S_1(\alpha)$.

\begin{lem}\label{lem:holder}
Let $\bxi =[\bxi_1,\ldots,\bxi_n]$ be a $K\times n$ random matrix with
real entries having a symmetric distribution (see the statement of \Cref{th:1}).
Let $\bzeta_i$ be defined as $\bzeta_i = \bxi_i \eta_i$, where
$\eta_i$ are iid random variables independent of $\bxi$ and satisfying
\begin{align}
\eta_i =
\begin{cases}
1,& \text{with probability }\ 1-\frac{\alpha\tau}{1+2\alpha\tau},\\
-1-\frac{1}{\alpha\tau},& \text{with probability }\ \frac{\alpha\tau}{1+2\alpha\tau}.
\end{cases}
\end{align}
Then, the expectation of the random variable $S$ can be bounded as follows:
\begin{align}
-(\nicefrac1\alpha)\,\bfE[S(\alpha)]  \le
\inf_{p}\bigg\{\int_{\calF} \ell_n(\bfF,\bfF^*)\,p(d\bfF) +
2\tau\KL(p\,||\pi_0)\bigg\},
\end{align}
where the inf is taken over all probability measures on $\calF$.
\end{lem}

\begin{lem}\label{lem:exp}
Let the random vectors $\bzeta_i$, $i\in[n]$ be as defined in \Cref{lem:holder}.
Then, we have
\begin{align}
\lim_{\alpha\to 0}\frac1{\alpha}\,\bfE[S_1(\alpha)\,|\,\bxi]\le
\sum_{i\in[n]} \tau\log \int_{\calF}e^{-(2/n\tau)\bxi_i^\top(\hat\bfF_i-\bfF_i)}\,\pi_n(d\bfF)
-\int_{\calF}\ell_n(\hat\bfF,\bfF)\pi_n(d\bfF).
\end{align}
\end{lem}

Applying these two lemmas, we get
\begin{align}
\bfE[\ell_n(\hat\bfF,\bfF^*)]
		& = \bfE[\ell_n(\hat\bfF,\bfF^*,\bzeta)]
		= \lim_{\alpha\to 0} \frac{\bfE(\bfE[S_1(\alpha)|\bxi])-\bfE[S(\alpha)]}{\alpha}\\
		&\le \inf_{p}\bigg\{\int_{\calF} \ell_n(\bfF,\bfF^*)\,p(d\bfF) + 2\tau\KL(p\,||\pi_0)\bigg\}
				-\int_{\calF}\ell_n(\hat\bfF,\bfF)\pi_n(d\bfF)\\
		&\qquad \qquad \qquad \qquad
		+ \sum_{i\in[n]} \tau\bfE\bigg\{\log \int_{\calF}e^{-(2/n\tau)\bxi_i^\top(\hat\bfF_i-\bfF_i)}
		 \,\pi_n(d\bfF)\bigg\}.\label{eq:3}
\end{align}
Then, for every $\tau\ge (2K/n)(B_\xi L)$, we have
\begin{align}
e^{-(2/n\tau)\bxi_i^\top(\hat\bfF_i-\bfF_i)}
		&\le 1 -\frac{2\bxi_i^\top(\hat\bfF_i-\bfF_i)}{n\tau}+
		\frac{3(\bxi_i^\top(\hat\bfF_i-\bfF_i))^2}{(n\tau)^2} \\
		&\le 1 -\frac{2\bxi_i^\top(\hat\bfF_i-\bfF_i)}{n\tau}+
		\frac{3KB_{\xi}^2\|\hat\bfF_i-\bfF_i\|_	2^2}{(n\tau)^2}.
\end{align}
This implies that
\begin{align}
\int_{\calF}e^{-(2/n\tau)\bxi_i^\top(\hat\bfF_i-\bfF_i)}
		 \,\pi_n(d\bfF)&\le 1 +
		\frac{3KB_{\xi}^2}{(n\tau)^2}\,\int_{\calF}\|\hat\bfF_i-\bfF_i\|_2^2\pi_n(d\bfF).
\end{align}
Combining the last display with \eqref{eq:3} and using the inequality
$\log(1+x)\le x$, we arrive at
\begin{align}
\bfE[\ell_n(\hat\bfF,\bfF^*)]
		\le \inf_{p}\bigg\{\int_{\calF} \ell_n(\bfF,\bfF^*)\,p(d\bfF) &+ 2\tau\KL(p\,||\pi_0)\bigg\}\\
		&-\Big(1-\frac{3KB_{\xi}^2}{n\tau}\Big)\bfE\bigg\{
		\int_{\calF}\ell_n(\hat\bfF,\bfF)\pi_n(d\bfF)\bigg\}.
\end{align}
This completes the proof of the theorem.
\end{proofof}

\begin{proofof}[Proof of \Cref{th:2}]
The proof follows the same arguments as those used in the proof of
\Cref{th:1}. That is why, we will skip some technical details. The main difference
is in the definition of the matrix $\bzeta$ and the subsequent computations related
to the evaluation of the term $S_2(\alpha)$. Thus, for any random matrix $\bzeta$ such
that  $\bfE[\bzeta|\bfY] = \mathbf 0$ and for
$\ell_n(\bfF,\bfF^*,\bzeta) = \ell_n(\bfF,\bfF^*) + \frac{2}{n}\langle\bzeta,\bfF-\bfF^*\rangle$,
we have
\begin{align}
\bfE[\ell_n(\hat\bfF,\bfF^*)] = \bfE[\ell_n(\hat\bfF,\bfF^*,\bzeta)]
		& = \lim_{\alpha\to 0}\frac{\bfE[S_1(\alpha)]-\bfE[S(\alpha)]}{\alpha},
\end{align}
where $S$ and $S_1$ are the same as in the proof of \Cref{th:1}. We instantiate the matrix
$\bzeta$ as follows: $\bzeta = \bfSigma^{1/2}\bar\bzeta$ where the entries of
$\bar\bzeta$ are given by $\bar\zeta_{j,i} = \bar\xi_{j,i} \eta_{j,i}$, with
$\eta_{j,i}$ being iid random variables independent of $\bxi$ and satisfying
\begin{align}
\eta_{j,i} =
\begin{cases}
1,& \text{with probability }\ 1-\frac{\alpha\tau}{1+2\alpha\tau},\\
-1-\frac{1}{\alpha\tau},& \text{with probability }\ \frac{\alpha\tau}{1+2\alpha\tau}.
\end{cases}
\end{align}
One easily checks that the resulting vector
$\bar\bxi_{j\bullet}+2\alpha\tau\bar\bzeta_{j,\bullet}$ has the same distribution as the
vector $(1+2\alpha\tau)\bar\bxi_{j,\bullet}$, for every $j\in[K]$. Furthermore, for different
values of $j$, these vectors are independent. This implies that the matrix
$\bar\bxi+2\alpha\tau\bar\bzeta$ has the same distribution as the matrix $(1+2\alpha\tau)\bar\bxi$,
which is sufficient for getting the conclusion of \Cref{lem:holder}. That is
\begin{align}
-(\nicefrac1\alpha)\,\bfE[S(\alpha)]  \le
\inf_{p}\bigg\{\int_{\calF} \ell_n(\bfF,\bfF^*)\,p(d\bfF) +
2\tau\KL(p\,||\pi_0)\bigg\},
\end{align}
where the inf is taken over all probability measures on $\calF$. To bound the term
$S_1$, we use a result similar to that of \Cref{lem:exp}.

\begin{lem}\label{lem:exp1}
Let the random matrix $\bzeta$ be defined as above.  Set $\bfH(\bfF) = \bfSigma^{1/2}
(\hat\bfF-\bfF)$. Then, we have
\begin{align}
\lim_{\alpha\to 0}\frac1{\alpha}\,\bfE[S_1(\alpha)\,|\,\bxi]\le
\sum_{\substack{i\in[n]\\ j\in[K]}}
\tau\log \int_{\calF}e^{-(2/n\tau)\bar\xi_{j,i}\bfH_{j,i}(\bfF)}\,\pi_n(d\bfF)
-\int_{\calF}\ell_n(\hat\bfF,\bfF)\pi_n(d\bfF).
\end{align}
\end{lem}
Then, for every $\tau\ge (2/n)(\bar B_{\xi}\bar L)$, we have
\begin{align}
e^{-(2/n\tau)\bar\xi_{j,i}\bfH_{j,i}(\bfF)}
		&\le 1 -\frac{2\bar\xi_{j,i}\bfH_{j,i}(\bfF)}{n\tau}+
		\frac{3\bar\xi_{j,i}^2\bfH_{j,i}(\bfF)^2}{(n\tau)^2} \\
		&\le 1 -\frac{2\bar\xi_{j,i}\bfH_{j,i}(\bfF)}{n\tau}+
		\frac{3\bar B_{\xi}^2\bfH_{j,i}^2(\bfF)}{(n\tau)^2}.
\end{align}
Using the fact that $\int_{\calF}\bfH(\bfF)\,\pi_n(d\bfF)=0$, the last display
implies that
\begin{align}
\int_{\calF}e^{-(2/n\tau)\bar\xi_{j,i}\bfH_{j,i}(\bfF)}
		 \,\pi_n(d\bfF)&\le 1 +
		\frac{3\bar B_{\xi}^2}{(n\tau)^2}\,\int_{\calF}\bfH_{j,i}^2(\bfF)\pi_n(d\bfF).
\end{align}
Combining the last display with \Cref{lem:exp1} and using the inequality
$\log(1+x)\le x$, we arrive at
\begin{align}
\lim_{\alpha\to 0}\frac1{\alpha}\,\bfE[S_1(\alpha)\,|\,\bxi]
		&\le \sum_{i,j} \frac{3\bar B_{\xi}^2}{n^2\tau}\,\int_{\calF}\bfH_{j,i}^2(\bfF)\pi_n(d\bfF)
			-\int_{\calF}\ell_n(\hat\bfF,\bfF)\pi_n(d\bfF)\\
		& = \frac{3\bar B_{\xi}^2}{n^2\tau}\,\int_{\calF}\|\bfSigma^{1/2}(\hat\bfF-\bfF)\|_F^2
				\pi_n(d\bfF)-\int_{\calF}\ell_n(\hat\bfF,\bfF)\pi_n(d\bfF)\\
		&\le \bigg(\frac{3\bar B_{\xi}^2\|\bfSigma\|}{n\tau}-1\bigg)\int_{\calF}\ell_n(\hat\bfF,\bfF)
		\pi_n(d\bfF).
\end{align}
This completes the proof of the theorem.
\end{proofof}

\begin{proofof}[Proof of \Cref{th:3}]
We outline here only the main steps of the proof, without going too much into the details.
One can extend the Stein lemma from the Gaussian distribution to that of $\bar\xi_{j,i}$,
provided the conditions of \Cref{th:3} are satisfied (see Lemma 1 in \citep{DT08} for a
similar result). The resulting claim is that the random variable
\begin{align}\label{eq:4}
\hat r:= \ell_n(\hat\bfF,\bfY) - \tr(\bfSigma) + \frac{2}{n}\sum_{i=1}^{n}\sum_{j=1}^K
g_{\bar\xi}(\bar\xi_{j,i})\partial_{\bar\xi_{j,i}} (\bfSigma^{1/2}\hat\bfF)_{j,i}
\end{align}
satisfies $\bfE[\hat r] = \bfE[\ell_n(\hat\bfF,\bfF^*)]$. On the one hand, using
Varadhan-Donsker's variational formula, we get
\begin{align}
\bfE[\ell_n(\hat\bfF,\bfY)]
		& \le \bfE[\ell_n(\hat\bfF,\bfY) + 2\tau \KL(\pi_n\,\|\,\pi_0)] \\
		& = \bfE\Big[\int_{\calF}\ell_n(\bfF,\bfY)\,\pi_n(d\bfF) + 2\tau \KL(\pi_n\,\|\,\pi_0)\Big]
		- \int_{\calF} \bfE\big[\ell_n(\bfF,\hat\bfF)\,\pi_n(d\bfF)\big]\\
		& \le \bfE\Big[\inf_p\Big(\int_{\calF}\ell_n(\bfF,\bfY)\,p(d\bfF) + 2\tau \KL(p\,\|\,\pi_0)
			\Big)\Big] - \int_{\calF} \bfE\big[\ell_n(\bfF,\hat\bfF)\,\pi_n(d\bfF)\big]\\
		& \le \inf_p\bfE\Big[\Big(\int_{\calF}\ell_n(\bfF,\bfY)\,p(d\bfF) + 2\tau \KL(p\,\|\,\pi_0)
			\Big)\Big] - \int_{\calF} \bfE\big[\ell_n(\bfF,\hat\bfF)\,\pi_n(d\bfF)\big]\\			
		& =\inf_p\Big(\ell_n(\bfF,\bfF^*)+\tr(\bfSigma) + 2\tau \KL(p\,\|\,\pi_0)\Big) -
			\int_{\calF} \bfE\big[\ell_n(\bfF,\hat\bfF)\,\pi_n(d\bfF)\big].\label{eq:5}
\end{align}
On the other hand, computing the partial derivative $\partial_{Y_{j,i}}
(\bfSigma^{1/2}\hat\bfF)_{j,i}$, we get
\begin{align}
\partial_{\bar\xi_{j,i}} (\bfSigma^{1/2}\hat\bfF)_{j,i}
		& = \be_j^\top\bfSigma^{1/2} (\partial_{\bar\bxi_{i}} \hat\bfF_i) \be_j\\
		& = \frac1{2n\tau} \be_j^\top\bfSigma^{1/2}\big(\int_{\calF}\bfF_i(\bfF_i-\bfY_i)^\top
				\pi_n(d\bfF)-\hat\bfF_i(\hat\bfF_i-\bfY_i)^\top\big)\bfSigma^{1/2}\be_j\\
		& = \frac1{2n\tau} \int_{\calF}\big\{\be_j^\top\bfSigma^{1/2}(\bfF-\hat\bfF)_{i}\big\}^2
		\,\pi_n(d\bfF).\label{eq:6}
\end{align}
From this relation, we infer that
\begin{align}
\sum_{j=1}^K
g_{\bar\xi}(\bar\xi_{j,i})\partial_{Y_{j,i}} (\bfSigma^{1/2}\hat\bfF)_{j,i}
		& = \frac1{2n\tau}\sum_{j=1}^K g_{\bar\xi}
		(\bar\xi_{j,i})\int_{\calF}\big\{\be_j^\top\bfSigma^{1/2}(\bfF-\hat\bfF)_{i}\big\}^2
		\pi_n(d\bfF)\\
		& \le \frac{G_{\bar\xi}}{2n\tau}\int_{\calF}\|\bfSigma^{1/2}
				(\bfF_i-\hat\bfF_i)\|_2^2\,\pi_n(d\bfF)\\
		&\le \frac{\|\bfSigma\| G_{\bar\xi}}{2n\tau} \int_{\calF}\|\bfF_i-\hat\bfF_i\|_2^2
				\,\pi_n(d\bfF).\label{eq:7}
\end{align}
Combining \eqref{eq:4}-\eqref{eq:7}, we arrive at
\begin{align}
\bfE[\ell_n(\hat\bfF,\bfF^*)]
		& = \bfE[\hat r] \\
		&\le \inf_p\Big(\ell_n(\bfF,\bfF^*) + 2\tau \KL(p\,\|\,\pi_0)\Big) -
			\int_{\calF} \bfE\big[\ell_n(\bfF,\hat\bfF)\,\pi_n(d\bfF)\big]\\
		&\qquad +
			\frac{\|\bfSigma\|G_{\bar\xi}}{n^2\tau} \sum_{i=1}^n \int_{\calF}
			\bfE[\|\hat\bfF_i-\bfF_i\|_2^2\,\pi_n(d\bfF)]\\
		&=\inf_p\Big(\ell_n(\bfF,\bfF^*) + 2\tau \KL(p\,\|\,\pi_0)\Big) -
			\bigg(1-\frac{\|\bfSigma\|G_{\bar\xi}}{n\tau}\bigg)\int_{\calF}
			\bfE\big[\ell_n(\bfF,\hat\bfF)\,\pi_n(d\bfF)\big].
\end{align}
This completes the proof.
\end{proofof}

\begin{proofof}[Proof of \Cref{prop:1}] Without loss of generality, we assume that $a\ge 0$.
We have, for every $x\in[a,b]$,
\begin{align}
m_{\bar\xi}(x) &= \int_x^b y p_{\bar\xi}(y)\,dy \le p_{\bar\xi}(x) \int_x^b y\,dy\le
(\nicefrac{b^2}{2}) p_{\bar\xi}(x).
\end{align}
Similarly, for  every $x\in[-b,0]$, we have $x\le a$ and, therefore,
\begin{align}
m_{\bar\xi}(x) &= -\int_{-b}^x y p_{\bar\xi}(y)\,dy \le p_{\bar\xi}(x) \int_{-b}^x (-y)\,dy\le (\nicefrac{b^2}{2}) p_{\bar\xi}(x).
\end{align}
Finally, for $x\in[0,a]$, we have
\begin{align}
m_{\bar\xi}(x) &= \int_{-b}^x (-y) p_{\bar\xi}(y)\,dy \le \int_{-b}^0 (-y) p_{\bar\xi}(y)\,dy
\le (\nicefrac{b^2}{2}) p_{\bar\xi}(0) \le (\nicefrac{b^2}{2}) p_{\bar\xi}(x)
\end{align}
and the claim of the lemma follows.
\end{proofof}

\acks{This work was partially supported
by the grant Investissements d'Avenir (ANR-11-IDEX-0003/Labex Ecodec/ANR-11-LABX-0047) and the chair
``LCL/GENES/Fondation du risque, Nouveaux enjeux pour nouvelles données''.}

\bibliography{bibDGP16}

\begin{thebibliography}{43}
\providecommand{\natexlab}[1]{#1}
\providecommand{\url}[1]{\texttt{#1}}
\expandafter\ifx\csname urlstyle\endcsname\relax
  \providecommand{\doi}[1]{doi: #1}\else
  \providecommand{\doi}{doi: \begingroup \urlstyle{rm}\Url}\fi

\bibitem[Alquier and Guedj(2017)]{Alquier2017}
P.~Alquier and B.~Guedj.
\newblock An oracle inequality for quasi-bayesian nonnegative matrix
  factorization.
\newblock \emph{Mathematical Methods of Statistics}, 26\penalty0 (1):\penalty0
  55--67, Jan 2017.
\newblock ISSN 1934-8045.
\newblock \doi{10.3103/S1066530717010045}.
\newblock URL \url{https://doi.org/10.3103/S1066530717010045}.

\bibitem[Alquier(2013)]{AlquierALT}
Pierre Alquier.
\newblock Bayesian methods for low-rank matrix estimation: Short survey and
  theoretical study.
\newblock In Sanjay Jain, R{\'e}mi Munos, Frank Stephan, and Thomas Zeugmann,
  editors, \emph{Algorithmic Learning Theory}, pages 309--323, Berlin,
  Heidelberg, 2013. Springer Berlin Heidelberg.

\bibitem[Alquier and Lounici(2011)]{Alquier3}
Pierre Alquier and Karim Lounici.
\newblock P{AC}-{B}ayesian bounds for sparse regression estimation with
  exponential weights.
\newblock \emph{Electron. J. Stat.}, 5:\penalty0 127--145, 2011.

\bibitem[Audibert(2009)]{Audibert}
Jean-Yves Audibert.
\newblock Fast learning rates in statistical inference through aggregation.
\newblock \emph{Ann. Statist.}, 37\penalty0 (4):\penalty0 1591--1646, 2009.

\bibitem[Bellec(2018)]{bellec2018}
Pierre~C. Bellec.
\newblock Optimal bounds for aggregation of affine estimators.
\newblock \emph{Ann. Statist.}, 46\penalty0 (1):\penalty0 30--59, 02 2018.
\newblock \doi{10.1214/17-AOS1540}.

\bibitem[Bouwmans et~al.(2016)Bouwmans, Aybat, and
  Zahzah]{bouwmans:hal-01373013}
Thierry Bouwmans, Necdet~Serhat Aybat, and El-Hadi Zahzah.
\newblock \emph{{Handbook on ''Robust Low-Rank and Sparse Matrix Decomposition:
  Applications in Image and Video Processing''}}.
\newblock {CRC Press, Taylor and Francis Group, }, May 2016.
\newblock URL \url{https://hal.archives-ouvertes.fr/hal-01373013}.

\bibitem[Bunea et~al.(2011{\natexlab{a}})Bunea, She, and Wegkamp]{BYW}
Florentina Bunea, Yiyuan She, and Marten~H. Wegkamp.
\newblock Optimal selection of reduced rank estimators of high-dimensional
  matrices.
\newblock \emph{Ann. Statist.}, 39\penalty0 (2):\penalty0 1282--1309,
  2011{\natexlab{a}}.

\bibitem[Bunea et~al.(2011{\natexlab{b}})Bunea, She, and Wegkamp]{bunea2011}
Florentina Bunea, Yiyuan She, and Marten~H. Wegkamp.
\newblock Optimal selection of reduced rank estimators of high-dimensional
  matrices.
\newblock \emph{Ann. Statist.}, 39\penalty0 (2):\penalty0 1282--1309, 04
  2011{\natexlab{b}}.
\newblock \doi{10.1214/11-AOS876}.
\newblock URL \url{https://doi.org/10.1214/11-AOS876}.

\bibitem[Bunea et~al.(2012)Bunea, She, and Wegkamp]{bunea2012}
Florentina Bunea, Yiyuan She, and Marten~H. Wegkamp.
\newblock Joint variable and rank selection for parsimonious estimation of
  high-dimensional matrices.
\newblock \emph{Ann. Statist.}, 40\penalty0 (5):\penalty0 2359--2388, 10 2012.
\newblock \doi{10.1214/12-AOS1039}.
\newblock URL \url{https://doi.org/10.1214/12-AOS1039}.

\bibitem[Cand{\`e}s and Plan(2011)]{CP11}
Emmanuel~J. Cand{\`e}s and Yaniv Plan.
\newblock Tight oracle inequalities for low-rank matrix recovery from a minimal
  number of noisy random measurements.
\newblock \emph{IEEE Trans. Inform. Theory}, 57\penalty0 (4):\penalty0
  2342--2359, 2011.

\bibitem[Cand{\`e}s and Tao(2010)]{CT10}
Emmanuel~J. Cand{\`e}s and Terence Tao.
\newblock The power of convex relaxation: near-optimal matrix completion.
\newblock \emph{IEEE Trans. Inform. Theory}, 56\penalty0 (5):\penalty0
  2053--2080, 2010.

\bibitem[Catoni(2007)]{Catoni1}
Olivier Catoni.
\newblock \emph{Pac-{B}ayesian supervised classification: the thermodynamics of
  statistical learning}.
\newblock Lecture Notes--Monograph Series, 56. Institute of Mathematical
  Statistics, Beachwood, OH, 2007.

\bibitem[Dai et~al.(2014)Dai, Rigollet, Xia, and Zhang]{DRXZ}
Dong Dai, Philippe Rigollet, Lucy Xia, and Tong Zhang.
\newblock Aggregation of affine estimators.
\newblock \emph{Electron. J. Stat.}, 8\penalty0 (1):\penalty0 302--327, 2014.

\bibitem[Dalalyan and Tsybakov(2012{\natexlab{a}})]{DT12b}
A.~S. Dalalyan and A.~B. Tsybakov.
\newblock Sparse regression learning by aggregation and {L}angevin
  {M}onte-{C}arlo.
\newblock \emph{J. Comput. System Sci.}, 78\penalty0 (5):\penalty0 1423--1443,
  2012{\natexlab{a}}.

\bibitem[Dalalyan(2017)]{Dal14}
Arnak~S. Dalalyan.
\newblock Theoretical guarantees for approximate sampling from a smooth and
  log-concave density.
\newblock to appear in JRSS B , arXiv:1412.7392, December 2017.

\bibitem[Dalalyan and Karagulyan(2017)]{DalKar17}
Arnak~S. Dalalyan and Avetik Karagulyan.
\newblock User-friendly guarantees for the langevin monte carlo with inaccurate
  gradient.
\newblock submitted 1710.00095, arXiv, October 2017.
\newblock URL \url{https://arxiv.org/abs/1710.00095}.

\bibitem[Dalalyan and Salmon(2012)]{DS12}
Arnak~S. Dalalyan and Joseph Salmon.
\newblock Sharp oracle inequalities for aggregation of affine estimators.
\newblock \emph{Ann. Statist.}, 40\penalty0 (4):\penalty0 2327--2355, 2012.

\bibitem[Dalalyan and Tsybakov(2007)]{DT07}
Arnak~S. Dalalyan and Alexandre~B. Tsybakov.
\newblock Aggregation by exponential weighting and sharp oracle inequalities.
\newblock In \emph{Learning theory}, volume 4539 of \emph{Lecture Notes in
  Comput. Sci.}, pages 97--111. Springer, Berlin, 2007.

\bibitem[Dalalyan and Tsybakov(2008)]{DT08}
Arnak~S. Dalalyan and Alexandre~B. Tsybakov.
\newblock Aggregation by exponential weighting, sharp pac-bayesian bounds and
  sparsity.
\newblock \emph{Machine Learning}, 72\penalty0 (1-2):\penalty0 39--61, 2008.

\bibitem[Dalalyan and Tsybakov(2012{\natexlab{b}})]{DT12a}
Arnak~S. Dalalyan and Alexandre~B. Tsybakov.
\newblock Mirror averaging with sparsity priors.
\newblock \emph{Bernoulli}, 18\penalty0 (3):\penalty0 914--944,
  2012{\natexlab{b}}.

\bibitem[Dalalyan et~al.(2016)Dalalyan, Grappin, and Paris]{DGP16}
Arnak~S. Dalalyan, Edwin Grappin, and Quentin Paris.
\newblock On the exponentially weighted aggregate with the laplace prior.
\newblock to appear in the Annals of Statistics 1611.08483, arXiv, November
  2016.
\newblock URL \url{https://arxiv.org/abs/1611.08483}.

\bibitem[Durmus and Moulines(2016)]{Durmus16}
A.~Durmus and E.~Moulines.
\newblock {Sampling from strongly log-concave distributions with the Unadjusted
  Langevin Algorithm}.
\newblock Technical Report , arXiv:1605.01559, May 2016.

\bibitem[Ga{\"{\i}}ffas and Lecu{\'e}(2011)]{GL11}
St{\'e}phane Ga{\"{\i}}ffas and Guillaume Lecu{\'e}.
\newblock Sharp oracle inequalities for high-dimensional matrix prediction.
\newblock \emph{IEEE Trans. Inform. Theory}, 57\penalty0 (10):\penalty0
  6942--6957, 2011.

\bibitem[George(1986{\natexlab{a}})]{George86a}
E.~I. George.
\newblock Minimax multiple shrinkage estimation.
\newblock \emph{Ann. Statist.}, 14\penalty0 (1):\penalty0 188--205,
  1986{\natexlab{a}}.

\bibitem[George(1986{\natexlab{b}})]{George86b}
E.~I. George.
\newblock Combining minimax shrinkage estimators.
\newblock \emph{J. Amer. Statist. Assoc.}, 81\penalty0 (394):\penalty0
  437--445, 1986{\natexlab{b}}.

\bibitem[Izenman(1975)]{IZENMAN1975248}
Alan~Julian Izenman.
\newblock Reduced-rank regression for the multivariate linear model.
\newblock \emph{Journal of Multivariate Analysis}, 5\penalty0 (2):\penalty0 248
  -- 264, 1975.
\newblock ISSN 0047-259X.
\newblock \doi{https://doi.org/10.1016/0047-259X(75)90042-1}.
\newblock URL
  \url{http://www.sciencedirect.com/science/article/pii/0047259X75900421}.

\bibitem[Juditsky et~al.(2008)Juditsky, Rigollet, and Tsybakov]{JRT}
A.~Juditsky, P.~Rigollet, and A.~B. Tsybakov.
\newblock Learning by mirror averaging.
\newblock \emph{Ann. Statist.}, 36\penalty0 (5):\penalty0 2183--2206, 2008.

\bibitem[Klopp(2014)]{Klopp14}
Olga Klopp.
\newblock Noisy low-rank matrix completion with general sampling distribution.
\newblock \emph{Bernoulli}, 20\penalty0 (1):\penalty0 282--303, 2014.

\bibitem[Koltchinskii et~al.(2011)Koltchinskii, Lounici, and Tsybakov]{KLT11}
Vladimir Koltchinskii, Karim Lounici, and Alexandre~B. Tsybakov.
\newblock Nuclear-norm penalization and optimal rates for noisy low-rank matrix
  completion.
\newblock \emph{The Annals of Statistics}, 39\penalty0 (5):\penalty0
  2302--2329, 2011.

\bibitem[Leung(2004)]{Leung}
G.~Leung.
\newblock \emph{Information Theory and Mixing Least Squares Regression}.
\newblock PhD thesis, Yale University, 2004.

\bibitem[Leung and Barron(2006)]{Leung2006}
Gilbert Leung and Andrew~R. Barron.
\newblock Information theory and mixing least-squares regressions.
\newblock \emph{IEEE Trans. Inform. Theory}, 52\penalty0 (8):\penalty0
  3396--3410, 2006.

\bibitem[Mai and Alquier(2015{\natexlab{a}})]{Mai15}
The~Tien Mai and Pierre Alquier.
\newblock A {B}ayesian approach for noisy matrix completion: optimal rate under
  general sampling distribution.
\newblock \emph{Electron. J. Stat.}, 9\penalty0 (1):\penalty0 823--841,
  2015{\natexlab{a}}.

\bibitem[Mai and Alquier(2015{\natexlab{b}})]{mai2015}
The~Tien Mai and Pierre Alquier.
\newblock A bayesian approach for noisy matrix completion: Optimal rate under
  general sampling distribution.
\newblock \emph{Electron. J. Statist.}, 9\penalty0 (1):\penalty0 823--841,
  2015{\natexlab{b}}.
\newblock \doi{10.1214/15-EJS1020}.
\newblock URL \url{https://doi.org/10.1214/15-EJS1020}.

\bibitem[Negahban and Wainwright(2011)]{NW11}
Sahand Negahban and Martin~J. Wainwright.
\newblock Estimation of (near) low-rank matrices with noise and
  high-dimensional scaling.
\newblock \emph{Ann. Statist.}, 39\penalty0 (2):\penalty0 1069--1097, 2011.

\bibitem[Negahban and Wainwright(2012)]{NW12}
Sahand Negahban and Martin~J. Wainwright.
\newblock Restricted strong convexity and weighted matrix completion: optimal
  bounds with noise.
\newblock \emph{J. Mach. Learn. Res.}, 13:\penalty0 1665--1697, 2012.

\bibitem[Rai et~al.(2014)Rai, Wang, Guo, Chen, Dunson, and Carin]{RaiICML}
Piyush Rai, Yingjian Wang, Shengbo Guo, Gary Chen, David Dunson, and Lawrence
  Carin.
\newblock Scalable bayesian low-rank decomposition of incomplete multiway
  tensors.
\newblock In Eric~P. Xing and Tony Jebara, editors, \emph{Proceedings of the
  31st International Conference on Machine Learning}, volume~32 of
  \emph{Proceedings of Machine Learning Research}, pages 1800--1808, Bejing,
  China, 22--24 Jun 2014. PMLR.

\bibitem[Rigollet and Tsybakov(2011)]{RT11}
Philippe Rigollet and Alexandre Tsybakov.
\newblock Exponential screening and optimal rates of sparse estimation.
\newblock \emph{Ann. Statist.}, 39\penalty0 (2):\penalty0 731--771, 2011.

\bibitem[Rigollet and Tsybakov(2012)]{RT12}
Philippe Rigollet and Alexandre~B. Tsybakov.
\newblock Sparse estimation by exponential weighting.
\newblock \emph{Statist. Sci.}, 27\penalty0 (4):\penalty0 558--575, 2012.

\bibitem[Rohde and Tsybakov(2011)]{Rohde11}
Angelika Rohde and Alexandre~B. Tsybakov.
\newblock Estimation of high-dimensional low-rank matrices.
\newblock \emph{Ann. Statist.}, 39\penalty0 (2):\penalty0 887--930, 2011.

\bibitem[Srebro and Shraibman(2005)]{Srebro2005}
Nathan Srebro and Adi Shraibman.
\newblock Rank, trace-norm and max-norm.
\newblock In Peter Auer and Ron Meir, editors, \emph{18th Annual Conference on
  Learning Theory, COLT 2005. Proceedings}, pages 545--560, 2005.

\bibitem[Yang et~al.(2017)Yang, Fang, Duan, Li, and Zeng]{Yang2017}
Linxiao Yang, Jun Fang, Huiping Duan, Hongbin Li, and Bing Zeng.
\newblock Fast low-rank bayesian matrix completion with hierarchical gaussian
  prior models.
\newblock \emph{CoRR}, abs/1708.02455, 2017.
\newblock URL \url{http://arxiv.org/abs/1708.02455}.

\bibitem[Yang(2000{\natexlab{a}})]{Yang00}
Y.~Yang.
\newblock Combining different procedures for adaptive regression.
\newblock \emph{J. Multivariate Anal.}, 74\penalty0 (1):\penalty0 135--161,
  2000{\natexlab{a}}.

\bibitem[Yang(2000{\natexlab{b}})]{Yang00a}
Y.~Yang.
\newblock Adaptive estimation in pattern recognition by combining different
  procedures.
\newblock \emph{Statist. Sinica}, 10\penalty0 (4):\penalty0 1069--1089,
  2000{\natexlab{b}}.

\end{thebibliography}

\newpage

\appendix

\section{Proofs of technical lemmas}\label{sec:lemmas}

\begin{proofof}[Proof of \Cref{lem:holder}]
Using the fact that $\pi_n(d\bfF)\propto\exp\big\{-(\nicefrac1{2\tau})\ell_n(\bfF,\bfY)\big\}\pi_0(d\bfF)$
and that $\ell_n(\bfF,\bfY) = \ell_n(\bfF,\bfF^*) + (\nicefrac{2}{n})\langle \bxi,\bfF^*-\bfF\rangle + (\nicefrac{1}{n})\|\bxi\|_F^2$, we arrive at
\begin{align}\label{L1:1}
-S(\alpha) & = \log \int_{\calF} e^{-(\nicefrac1{2\tau})\ell_n(\bfF,\bfF^*)-
		(\nicefrac{1}{n\tau})\langle \bxi,\bfF^*-\bfF\rangle}\,\pi_0(d\bfF)\\
		&\qquad -
		\log \int_{\calF} e^{-(\alpha+\nicefrac1{2\tau})\ell_n(\bfF,\bfF^*)-
		(\nicefrac{1}{n\tau})\langle \bxi + 2\alpha\tau\bzeta,\bfF^*-\bfF\rangle}\,\pi_0(d\bfF)
\end{align}
One easily checks that the random matrix $\bxi+2\alpha\tau\bzeta$ has the same distribution as the
matrix $(1+2\alpha\tau)\bxi$ and, therefore,
\begin{align}\label{L1:2}
-\bfE[S(\alpha)] & = \bfE\Big[\log \int_{\calF} e^{-(\nicefrac1{2\tau})\ell_n(\bfF,\bfF^*)-
		(\nicefrac{1}{n\tau})\langle \bxi,\bfF^*-\bfF\rangle}\,\pi_0(d\bfF)\Big]\\
		&\qquad -\bfE\Big[
		\log \int_{\calF} e^{-(2\alpha\tau +1)\big(\nicefrac1{2\tau}\ell_n(\bfF,\bfF^*)+
		(\nicefrac{1}{n\tau})\langle \bxi,\bfF^*-\bfF\rangle\big)}\,\pi_0(d\bfF)\Big].
\end{align}
Applying the H\"older inequality
$\int_{\calF} G\,d\pi_0\le \big(\int_{\calF} G^{2\alpha\tau+1}\,d\pi_0\big)^{1/(2\alpha\tau+1)}$ to
the first expectation of the right hand side, we get
\begin{align}\label{L1:3}
-\bfE[S(\alpha)] & \le -\frac{2\alpha\tau}{2\alpha\tau+1}\bfE\Big[
		\log \int_{\calF} e^{-(2\alpha\tau +1)\big(\nicefrac1{2\tau}\ell_n(\bfF,\bfF^*)+
		(\nicefrac{1}{n\tau})\langle \bxi,\bfF^*-\bfF\rangle\big)}\,\pi_0(d\bfF)\Big].
\end{align}
Donsker-Varadhan's variational inequality implies that
\begin{align}\label{L1:4}
-\frac1\alpha\,\bfE[S(\alpha)]
		& \le \bfE\Big[\inf_p\Big\{\int_{\calF} \big(\ell_n(\bfF,\bfF^*)+
				(\nicefrac{2}{n})\langle \bxi,\bfF^*-\bfF\rangle\big)\,p(d\bfF) +
				\frac{2\tau}{2\alpha\tau+1}\,\KL(p\,||\pi_0)\Big\}\Big]\\
		&\le  \inf_p\Big\{\int_{\calF} \bfE\big[\ell_n(\bfF,\bfF^*)+
				(\nicefrac{2}{n})\langle \bxi,\bfF^*-\bfF\rangle\big]\,p(d\bfF) +
				\frac{2\tau}{2\alpha\tau+1}\,\KL(p\,||\pi_0)\Big\}\\
		&\le  \inf_p\Big\{\int_{\calF} \bfE\big[\ell_n(\bfF,\bfF^*)\big]\,p(d\bfF) +
				\frac{2\tau}{2\alpha\tau+1}\,\KL(p\,||\pi_0)\Big\}	.
\end{align}
The desired result follows from the last display using the inequality
$2\alpha\tau +1 \ge 1$.
\end{proofof}

\begin{proofof}[Proof of \Cref{lem:exp}]
We have
\begin{align}
\ell_n(\hat\bfF,\bfF^*,\bzeta)-\ell_n(\bfF,\bfF^*,\bzeta)
		= \ell_n(\hat\bfF,\bfF^*)-\ell_n(\bfF,\bfF^*) + \frac2n \langle \bzeta,\hat\bfF- \bfF\rangle
\end{align}
which implies that,
\begin{align}
S_1(\alpha)
			&= \log \int_{\calF}e^{\alpha\big(\ell_n(\hat\bfF,\bfF^*)-\ell_n(\bfF,\bfF^*) +
					\frac2n \langle \bzeta,\hat\bfF- \bfF\rangle\big)}\,\pi_n(d\bfF).
\end{align}
Using the definition of the expectation, we get
\begin{align}
\Psi(\alpha)
			&:=\bfE[S_1(\alpha)\,|\,\bxi\,]\\
			&= \sum_{\bss\in \{0,1\}^n} \frac{(\alpha\tau)^{\|\bss\|_1}
			(1+\alpha\tau)^{n-\|\bss\|_1}}{(1+2\alpha\tau)^n}			
			\log \int_{\calF}e^{\Phi(\alpha,\bss,\bfF)}\,\pi_n(d\bfF),
\end{align}
where
\begin{align}
\Phi(\alpha,\bss,\bfF)
			&:= \alpha\big(\ell_n(\hat\bfF,\bfF^*)-\ell_n(\bfF,\bfF^*)\big) +
					\frac2n \sum_{i=1}^n \Big\{\alpha(1-s_i)-s_i\Big(\alpha+\frac{1}{\tau}\Big)\Big\}\langle \bxi_i,\hat\bfF_i- \bfF_i\rangle.
\end{align}
One easily checks that the function $\Psi(\alpha)$ is differentiable in $(0,\infty)$ and
that $\Psi(0) = 0$. Therefore,
\begin{align}
\lim_{\alpha\to0}\frac{\Psi(\alpha)}{\alpha}
		& = \lim_{\alpha\to0}\frac{\Psi(\alpha)-\Psi(0)}{\alpha} = \Psi'(0)\\
		& = \frac{d}{d\alpha}\bigg|_{\alpha=0}\sum_{\substack{\bss\in \{0,1\}^n\\ \|\bss\|_1\le 1}}
				\frac{(\alpha\tau)^{\|\bss\|_1}
			(1+\alpha\tau)^{n-\|\bss\|_1}}{(1+2\alpha\tau)^n}			
			\log \int_{\calF}e^{\Phi(\alpha,\bss,\bfF)}\,\pi_n(d\bfF).\label{eq:2}
\end{align}
Let us now compute the derivatives with respect to $\alpha$ of the terms of the last sum.
For the term corresponding to $\bss =\mathbf 0$, since ${\Phi(0,\mathbf 0,\bfF)}=0$,
we have
\begin{align}
\frac{d}{d\alpha}\bigg|_{\alpha=0}
				\bigg\{\frac{(1+\alpha\tau)^n}{(1+2\alpha\tau)^n}			
			\log \int_{\calF}e^{\Phi(\alpha,\mathbf 0,\bfF)}\,\pi_n(d\bfF)\bigg\}
			&= \frac{d}{d\alpha}\bigg|_{\alpha=0}\bigg\{
			\log \int_{\calF}e^{\Phi(\alpha,\mathbf 0,\bfF)}\,\pi_n(d\bfF)\bigg\}\\
			& = \frac{d}{d\alpha}\bigg|_{\alpha=0}\bigg\{
			\int_{\calF}e^{\Phi(\alpha,\mathbf 0,\bfF)}\,\pi_n(d\bfF)\bigg\}\\
			& =
			\int_{\calF}\frac{d}{d\alpha}\bigg|_{\alpha=0} \Phi(\alpha,\mathbf 0,\bfF)\,\pi_n(d\bfF).
\end{align}
 Using that $\Phi(\alpha,\mathbf 0,\bfF)$ is a linear function of $\alpha$, as well as
the fact that $\int \bfF\,\pi_n(d\bfF) = \hat\bfF$, we arrive at
\begin{align}
\frac{d}{d\alpha}\bigg|_{\alpha=0}
				\bigg\{\frac{(1+\alpha\tau)^n}{(1+2\alpha\tau)^n}			
			\log \int_{\calF}e^{\Phi(\alpha,\mathbf 0,\bfF)}\,\pi_n(d\bfF)\bigg\}
			&=  -\int_{\calF} \ell_n(\hat\bfF,\bfF)\,\pi_n(d\bfF).
\end{align}
Let us go back to \eqref{eq:2} and evaluate the terms corresponding to vectors $\bss$
such that $\|\bss\|_1=1$. This means that only one entry of $\bss$ is equal to one, all
the others being equal to zero. Hence, if we denote by $\be_i$ the $i$th element of the canonical
basis of $\RR^n$, we get
\begin{align}
\frac{d}{d\alpha}\bigg|_{\alpha=0}&\sum_{\substack{\bss\in \{0,1\}^n\\ \|\bss\|_1= 1}}
				\frac{(\alpha\tau)^{\|\bss\|_1}
			(1+\alpha\tau)^{n-\|\bss\|_1}}{(1+2\alpha\tau)^n}			
			\log \int_{\calF}e^{\Phi(\alpha,\bss,\bfF)}\,\pi_n(d\bfF)\\
			& =\sum_{i\in[n]} \frac{d}{d\alpha}\bigg|_{\alpha=0}
				\frac{(\alpha\tau)(1+\alpha\tau)^{n-1}}{(1+2\alpha\tau)^n}			
			\log \int_{\calF}e^{\Phi(\alpha,\be_i,\bfF)}\,\pi_n(d\bfF)\\
			& =\sum_{i\in[n]} \tau\log \int_{\calF}e^{\Phi(0,\be_i,\bfF)}\,\pi_n(d\bfF)	
				=\sum_{i\in[n]} \tau\log \int_{\calF}e^{-(2/n\tau)\bxi_i^\top(\hat\bfF_i-\bfF_i)}\,\pi_n(d\bfF).
\end{align}			
This completes the proof of the lemma.
\end{proofof}

\begin{proofof}[Proof of \Cref{lem:t3}]
For any bounded and measurable function $h:\RR^K\to\RR$, we have
\begin{align}
\int_{\calF}h(\bfF_1)\,\pi_0(\bfF)\,d\bfF
		& \stackrel{(a)}{=} \frac1{C_{\lambda}} \int_{\calF}
		\frac{h(\bfF_1)}{\text{det}(\lambda^2\bfI_K+\bfF\bfF^\top)^{(n+K+2)/2}} d\bfF\\
		& \stackrel{(b)}{=}  \frac{1}{C_1}
		\int_{\calF}\frac{h(\lambda\bfM_1)}{\text{det}(\bfI_K+ \bfM\bfM^\top)^{(n+K+2)/2}} d\bfM
\end{align}
where in (a) we have used the notation $C_\lambda = \int_{\calF}
\text{det}(\lambda^2\bfI_K+\bfF\bfF^\top)^{-(n+K+2)/2} d\bfF$, whereas in (b) we have
made the change of variable $\bfF = \lambda\bfM$. In order to compute the last integral,
we make another change of variable, $\bfM \rightsquigarrow \bar\bfM$, given by
$\bfM = [\bar\bfM_1 , (\bfI+\bar\bfM_1\bar\bfM_1^\top)^{1/2}\bar\bfM_{2:n}]$. This yields
\begin{align}
d\bfM = d\bfM_1 d\bfM_{2:n}
		& = d\bar\bfM_1 \text{det}(\bfI+\bar\bfM_1\bar\bfM_1^\top)^{(n-1)/2}d\bar\bfM_{2:n}\\
		& \stackrel{(c)}{=} (1+\|\bar\bfM_1\|_2^2)^{(n-1)/2}d\bar\bfM_1d\bar\bfM_{2:n}
\end{align}
and
\begin{align}		
\text{det}(\bfI+\bfM\bfM^\top)
		& = \text{det}\big(\bfI+\bfM_1\bfM_1^\top + \bfM_{2:n}\bfM_{2:n}^\top\big)\\
		& = \text{det}\big(\bfI+\bar\bfM_1\bar\bfM_1^\top + (\bfI+\bar\bfM_1\bar\bfM_1^\top)^{1/2}
			\bar\bfM_{2:n}\bar\bfM_{2:n}^\top(\bfI+\bar\bfM_1\bar\bfM_1^\top)^{1/2}\big)\\
		& = \text{det}\big((\bfI+\bar\bfM_1\bar\bfM_1^\top)^{1/2}(\bfI+
			\bar\bfM_{2:n}\bar\bfM_{2:n}^\top)(\bfI+\bar\bfM_1\bar\bfM_1^\top)^{1/2}\big)\\
		& = \text{det}(\bfI+\bar\bfM_1\bar\bfM_1^\top)\text{det}(\bfI+
			\bar\bfM_{2:n}\bar\bfM_{2:n}^\top)\\
		& \stackrel{(c')}{=} (1+\|\bar\bfM_1\|_2^2)\text{det}(\bfI+
			\bar\bfM_{2:n}\bar\bfM_{2:n}^\top),
\end{align}
where in (c) and (c$'$) we have used the fact that the matrix $\bfI+\bar\bfM_1\bar\bfM_1^\top$
has all its eigenvalues equal to 1 except the largest one, corresponding to the eigenvector
$\bar\bfM_1$, which is equal to $1+\|\bar\bfM_1\|_2^2$.
Using the same change of variable in $C_1$, and replacing $\bar\bfM_1$ by $\bx$ for convenience,
we get
\begin{align}
\int_{\calF}h(\bfF_1)\,\pi_0(\bfF)\,d\bfF
		& = \frac{\int_{\RR^K}h(\lambda\bx)(1+\|\bx\|_2^2)^{-(K+3)/2}\,d\bx}
		{\int_{\RR^K}(1+\|\bx\|_2^2)^{-(K+3)/2}\,d\bx}\\
		& = \frac{\int_{\RR^K}h(\lambda\by/\sqrt{3})(1+\|\by\|_2^2/3)^{-(K+3)/2}\,d\by}
		{3\int_{\RR^K}(1+\|\by\|_2^2/3)^{-(K+3)/2}\,d\by}.
\end{align}
In the last expression, we recognize the probability density function of the multivariate
$t_3$-distribution. Since the covariance matrix of a $K$-variate $t_\nu$ distribution
is equal to $\frac{\nu}{\nu-2}\bfI_K$, we infer that
\begin{align}
\int_{\calF}\|\bfF\|_F^2\,\pi_0(\bfF)\,d\bfF = nK\lambda^2.
\end{align}
This completes the proof of the lemma.
\end{proofof}

\begin{proofof}[Proof of \Cref{lem:KL}]
It holds that
\begin{align}
\KL(\bar p\|\pi_0)
		& = \int_{\calF} \log \Big(\frac{\pi_0(\bfF)}{\bar p(\bfF)}\Big)\,\pi_0(\bfF)\,d\bfF\\
		& = \int_{\calF} \log \Big(\frac{\pi_0(\bfF)}{\pi_0(\bfF-\bar\bfF)}\Big)\,\pi_0(\bfF)\,d\bfF.
\end{align}
To ease notation, we set $\bfA = (\lambda^2\bfI_K + \bfF\bfF^\top)^{-1/2}$ and
$\bfB = \lambda^2\bfI_K + (\bfF-\bar\bfF)(\bfF-\bar\bfF)^\top$. We have
\begin{align}\label{eq:9}
2\log \Big(\frac{\pi_0(\bfF)}{\pi_0(\bfF-\bar\bfF)}\Big)
		& = (n+K+2) \log\Big(\frac{\text{det}(\bfB)}{\text{det}(\bfA^{-2})}\Big)\\
		& = (n+K+2) \log\big(\text{det}(\bfA\bfB\bfA)\big)\\
		&	= (n+K+2) \sum_{j=1}^K\log s_j(\bfA\bfB\bfA),
\end{align}
where $s_j(\bfA\bfB\bfA)$ is the $j$th largest eigenvalue of the symmetric matrix
$\bfA\bfB\bfA$. Let $r$ be the rank of $\bar\bfF$. The first claim is that the matrix
$\bfA\bfB\bfA$ has at most $2r$ singular values different from one. Indeed, one can check
that
\begin{align}
\bfA\bfB\bfA -\bfI_K = \bfA\bar\bfF\bfF^*{}^\top\bfA - \bfA\bar\bfF\bfF^\top\bfA -
\bfA\bfF\bar\bfF{}^\top\bfA.
\end{align}
The matrix at the right hand side is at most of rank $2r$. This implies that
$\bfA\bfB\bfA -\bfI_K$ has at most $2r$ nonzero eigenvalues. Therefore, the number
of eigenvalues of $\bfA\bfB\bfA$ different from $1$ is not larger than $2r$, which
implies that the sum at the right hand side of \eqref{eq:9} has at most $2r$
nonzero entries.

Let $\bu_j$ be the unit eigenvector corresponding to the eigenvalue $s_j$. We know that, for
every $j\in[2r]$, $s_j = \bu_j^\top \bfA\bfB\bfA\bu_j$. Using \label{eq:8}, we get
\begin{align}
s_j &= 1+ \bu_j^\top(\bfA\bar\bfF\bfF^*{}^\top\bfA - \bfA\bar\bfF\bfF^\top\bfA -
				\bfA\bfF\bar\bfF{}^\top\bfA)\bu_j\\
		&= 1 - \|\bar\bfF{}^\top\bfA\bu_j\|_2^2+ \|(\bar\bfF-\bfF)^\top\bfA\bu_j\|_2^2\\
		&\le (1+\|\bar\bfF{}^\top\bfA\bu_j\|_2)^2.
\end{align}
Using the concavity of the function $\log(1+x^{1/2})$ over $(0,+\infty)$, we arrive at
\begin{align}\label{eq:10}
2\log \Big(\frac{\pi_0(\bfF)}{\pi_0(\bfF-\bar\bfF)}\Big)
		& = (n+K+2) \sum_{j=1}^{2r}\log s_j(\bfA\bfB\bfA)\\
		& \le 2(n+K+2)\sum_{j=1}^{2r}\log \big(1+\big\{\|\bar\bfF{}^\top\bfA\bu_j\|_2^2\big\}^{1/2}\big)\\
		& \le 4(n+K+2)r \log\bigg(1 + \bigg\{\frac1{2r}\sum_{j=1}^{2r}\|\bar\bfF{}^\top\bfA\bu_j\|_2^2
				\bigg\}^{1/2}\bigg) .
\end{align}
Finally, since $\bu_j$'s are orthonormal and $\bfA\preceq \lambda^{-1}\bfI_K$, we get
the claim of the lemma.
\end{proofof}

\section{Flaw in Corollary 1 of \citep{DT08}}\label{app:B}

As mentioned in the introduction, Corollary 1 in \citep{DT08}  relies heavily on
Lemma 3 of the same paper, that claims that
\begin{align}
x+ \log\bigg(1+\frac1{\alpha_0}(e^{-x\alpha_0}-1)\bigg) \le \frac{x^2\alpha_0}{2},\qquad
\forall x\in\RR,\ \forall \alpha_0>0.
\end{align}
Unfortunately, this inequality is not always true. In particular, the argument of
the logarithm is not always positive, which implies that the left hand side is not
always well defined. For instance, one can check that if $\alpha = 0.5$ and $x\ge 2$,
we have
\begin{align}
1+\frac1{\alpha_0}(e^{-x\alpha_0}-1)  =  2e^{-0.5x}-1 \le (2/e)-1 \le 0.
\end{align}

\end{document}